\documentclass[11pt]{article}
 \usepackage{benstyle}
 
\newcommand*{\pbk}{\vspace{1pc}\noindent} 
\newtheoremup{algorithm}[theorem]{Algorithm}
\usepackage{url}
\usepackage[font=small]{caption}

\title{
Optimal sampling rates for approximating analytic functions from pointwise samples
}
\author{Ben Adcock\\ Department of Mathematics \\ Simon Fraser University \\ Canada \\ \hspace{20pc} \\ \vspace{-2.25pc} \\ \and Rodrigo B. Platte \\ School of Mathematical and Statistical Sciences \\ Arizona State University \\ USA \and Alexei Shadrin \\ DAMTP, Centre for Mathematical Sciences \\ University of Cambridge \\ United Kingdom}
\begin{document}

\maketitle

\begin{abstract}
We consider the problem of approximating an analytic function on a compact interval from its values at $M+1$ distinct points.  When the points are equispaced, a recent result (the so-called \textit{impossibility theorem}) has shown that the best possible convergence rate of a stable method is root-exponential in $M$, and that any method with faster exponential convergence must also be exponentially ill-conditioned at a certain rate.  This result hinges on a classical theorem of Coppersmith \& Rivlin concerning the maximal behaviour of polynomials bounded on an equispaced grid.  In this paper, we first generalize this theorem to arbitrary point distributions.  We then present an extension of the impossibility theorem valid for general nonequispaced points, and apply it to the case of points that are equidistributed with respect to (modified) Jacobi weight functions.  This leads to a necessary sampling rate for stable approximation from such points.  We prove that this rate is also sufficient, and therefore exactly quantify (up to constants) the precise sampling rate for approximating analytic functions from such node distributions with stable methods.  Numerical results -- based on computing the maximal polynomial via a variant of the classical Remez algorithm -- confirm our main theorems.  
Finally, we discuss the implications of our results for polynomial least-squares approximations.  In particular, we theoretically confirm the well-known heuristic that stable least-squares approximation using polynomials of degree $N < M$ is possible only once $M$ is sufficiently large for there to be a subset of $N$ of the nodes that mimic the behaviour of the $N^{\rth}$ set of Chebyshev nodes.
\end{abstract}

\section{Introduction}
The concern of this paper is the approximation of an analytic function $f : [-1,1] \rightarrow \bbC$ from its values on an arbitrary set of $M+1$ points in $[-1,1]$.  It is well known that if such points follow a Chebyshev distribution then $f$ can be stably (up to a log factor in $M$) approximated by its polynomial interpolant with a convergence rate that is geometric in the parameter $M$.  Conversely, when the points are equispaced polynomial interpolants do not necessarily converge uniformly on $[-1,1]$ as $M \rightarrow \infty$; an effect known as Runge's phenomenon.  Such an approximation is also exponentially ill-conditioned in $M$, meaning that divergence is witnessed in finite precision arithmetic even when theoretical convergence is expected. 

Many different numerical methods have been proposed to \textit{overcome} Runge's phenomenon by replacing the polynomial interpolant by an alternative approximation (see \cite{BoydRunge,PlaKle,TrefPlatteIllCond} and references therein).  This raises the fundamental question: how successful can such approximations be?  For equispaced points, this question was answered recently in \cite{TrefPlatteIllCond}.  Therein it was proved that no stable method for approximating analytic functions from equispaced nodes can converge better than root-exponentially fast in the number of points, and moreover any method that converges exponentially fast must also be exponentially ill-conditioned.  

A well-known method for approximating analytic functions is polynomial least-squares fitting with a polynomial of degree $N < M$.  Although a classical approach, this technique has become increasingly popular in recent years as a technique for computing so-called polynomial chaos expansions with application to uncertainty quantification (see \cite{ChkifaEtAl,DavenportEtAlLeastSquares,MiglioratiEtAlFoCM,NarayanZhouCCP} and references therein), as well as in data assimilation in reduced-order modelling \cite{BinevEtAlDataAssim,DeVoreEtAlDataAssimBanach}.  A consequence of the impossibility theorem of \cite{TrefPlatteIllCond} is that polynomial least-squares is an optimal stable and convergent method for approximating one-dimensional analytic functions from equispaced data, provided the polynomial degree $N$ used in the least-squares fit scales like the square-root of the number of grid points $M+1$ \cite{Adcockl1Pointwise}.  Using similar ideas, it has also recently been shown that polynomial least-squares is also an optimal, stable method for extrapolating analytic functions \cite{DemanetTownsendExtrap}.

\subsection{Contributions}
The purpose of this paper is to investigate the limits of stability and accuracy for approximating analytic functions from arbitrary sets of points.  Of particular interest is the case of points whose behaviour lies between the two extremes of Chebyshev and equispaced grids.  Specifically, suppose a given set of points exhibits some clustering near the endpoints, but not the characteristic quadratic clustering of Chebyshev grids.  Generalizing that of \cite{TrefPlatteIllCond}, our main result precisely quantifies both the best achievable error decay rate for a stable approximation and the resulting ill-conditioning if one seeks faster convergence.  

This result follows from an extension of a classical theorem of Coppersmith \& Rivlin on the maximal behaviour of a polynomial of degree $N$ bounded on a grid of $M$ equispaced points \cite{copprivlinpolygrowth}.  In Theorem \ref{t:CoppRiv} we extend the lower bound proved in \cite{copprivlinpolygrowth} to arbitrary sets of points.  We next present an abstract impossibility result (Lemma \ref{l:absimpossibility}) valid for arbitrary sets of points.  To illustrate this result in a concrete setting, we then specialize to the case of nodes which are equidistributed with respect to modified Jacobi weight functions.  Such weight functions take the form
\bes{
\mu(x) = g(x) (1-x)^{\alpha} (1+x)^{\beta},
}
where $\alpha , \beta > -1$ and $c_1 \leq g(x) \leq c_2$ almost everywhere, and include equispaced ($\mu(x) = \frac12$) and Chebyshev ($\mu(x) = \frac{1}{\pi \sqrt{1-x^2}}$) weight functions as special cases.  In our main result, Theorem \ref{t:extendedimpossibility}, we prove an extended impossibility theorem for the corresponding nodes.  Generalizing \cite{TrefPlatteIllCond}, two important consequences of this theorem are as follows:
\begin{itemize}
\item[(i)] If $\gamma = \max \{ \alpha,\beta \} > -1/2$, any method that converges exponentially fast in $M$ with geometric rate, i.e.\ the error decays like $\ordu{\rho^{-M}}$ for some $\rho > 1$, must also be exponentially ill-conditioned in $M$ at a geometric rate.
\item[(ii)] The best possible convergence rate for a stable approximation is subgeometric with index $\nu = \frac{1}{2(\gamma+1)}$.  That is, the error is at best $\ordu{\rho^{-M^{\nu}}}$ for some $\rho > 1$ as $M \rightarrow \infty$.
\end{itemize}
We also give a full characterization of the trade-off between exponential ill-conditioning and exponential convergence at subgeometric rates lying strictly between $\nu = \frac{1}{2(\gamma+1)}$ and $\nu = 1$.

Although not a result about polynomials per se, this theorem is closely related to the behaviour of discrete least-squares fitting with polynomials of degree $N < M$.  Indeed, the quantity we estimate in Theorem \ref{t:CoppRiv} is equivalent (up to a factor of $\sqrt{M}$) to the infinity-norm condition number of such an approximation.  By using polynomial least-squares as our method, in Proposition \ref{p:least_squares_optimal} we show that the rate described in (ii) is not only necessary for stable recovery but also sufficient.  Specifically, when the polynomial degree is chosen as
\be{
\label{stable_scaling}
N \asymp M^{\nu},\qquad \nu = \frac{1}{2(\gamma+1)},
}
the polynomial least-squares approximation is stable and converges like $\ordu{\rho^{-M^{\nu}}}$ for all functions analytic in an appropriate complex region.  The fact that such a method is optimal in view of the generalized impossibility theorem goes some way towards justifying the popularity of discrete least-squares techniques.

Besides these results, in \S \ref{s:Remez} we also introduce an algorithm for computing the maximal polynomial for an arbitrary set of nodes.  This algorithm, which is based on a result of Sch\"onhage \cite{schonhagepolyinterp}, is a variant of the classical Remez algorithm for computing best uniform approximations (see, for example, \cite{PachonTrefethenRemez,powell}).  We use this algorithm to present numerical results in the paper confirming our various theoretical estimates.

Finally, let us note that one particular consequence of our results is a confirmation of a popular heuristic for polynomial least-squares approximation (see, for example, \cite{boyd2009divergence}).  Namely, the number of nonequispaced nodes required to stably recover a polynomial approximation of degree $N$ is of the same order as the number of nodes required for there to exist a subset of those nodes of size $N+1$ which mimics the distribution of the Chebyshev nodes $\{ \cos(n \pi/N) \}^{N}_{n=0}$.  In Proposition \ref{p:interlace} we show that the same sufficient condition for boundedness of the maximal polynomial also implies the existence of a subset of $N$ of the original $M+1$ nodes which interlace the Chebyshev nodes.  In particular, for nodes that are equidistributed according to a modified Jacobi weight function one has this interlacing property whenever the condition $M \asymp N^{2(\gamma+1)}$ holds, which is identical to the necessary and sufficient condition \R{stable_scaling} for stability of the least-squares approximation.

\section{Preliminaries}
Our focus in this paper is on functions defined on compact intervals, which we normalize to the unit interval $[-1,1]$.  Unless otherwise stated, all functions will be complex-valued, and in particular, polynomials may have complex coefficients.  Throughout the paper $-1 = x_0 < \ldots < x_M = 1$ will denote the nodes at which a function $f : [-1,1] \rightarrow \bbC$ is sampled.  We include both endpoints $x = \pm 1$ in this set for convenience.  All results we prove remain valid (with minor alterations) for the case when either or both endpoints is excluded.  

We require several further pieces of notation.  Where necessary throughout the paper $N \leq M$ will denote the degree of a polynomial.  We write $\| f \|_{\infty}$ for the uniform norm of a function $f \in C([-1,1])$ and $\| f \|_{M,\infty} = \max_{m=1,\ldots,M} | f(x_m) |$ for the discrete uniform semi-norm of $f$ on the grid of points.   We write $\ip{\cdot}{\cdot}$ for the Euclidean inner product on $L^2(-1,1)$ and $\nm{\cdot}_{2}$ for the Euclidean norm.  Correspondingly, we let $\ip{f}{g}_{M} = \frac{2}{M+1} \sum^{M}_{m=0} f(x_m) \overline{g(x_m)}$ and $\| f \|_{M,2} = \sqrt{\ip{f}{f}_M}$ be the discrete semi-inner product and semi-norm respectively.

We will say that a sequence $a_n$ converges to zero exponentially fast if $|a_n| = \ord{\rho^{-n^r}}$ for some $\rho > 1$ and $r > 0$.  If $r=1$ then we say the convergence is geometric, and if $0 < r <1$ or $r >1$ then it is subgeometric or supergeometric respectively.  When $r=1/2$ we also refer to this convergence as root-exponential.  Given two nonnegative sequences $a_n$ and $b_n$ we write $a_n \asymp b_n$ as $n \rightarrow \infty$ if there exist constants $c_1,c_2 >0$ such that $c_1 b_n \leq a_n \leq c_2 b_n$ for all large $n$.  Finally, we will on occasion use the notation $A \lesssim B$ to mean that there exists a constant $c>0$ independent of all relevant parameters such that $A \leq c B$.

\subsection{The impossibility theorem for equispaced points}
We first review the impossibility theorem of \cite{TrefPlatteIllCond}.  Let $\{ x_m \}^{M}_{m=0} = \{ -1 + 2 m/M \}^{M}_{m=0}$ be a grid of $M+1$ equispaced points in $[-1,1]$ and suppose that $F_{M} : C([-1,1]) \rightarrow C([-1,1])$ is a family of mappings such that $F_{M}(f)$ depends only on the values of $f$ on this grid.  We define the (absolute) condition numbers as
\be{
\label{condition_number}
\kappa(F_M) = \sup_{f \in C([-1,1])} \lim_{\delta \rightarrow 0^{+}} \sup_{\substack{h \in C([-1,1]) \\ 0 < \| h \|_{M,\infty} \leq \delta}} \frac{\| F_M(f+h) - F_M(f) \|_{\infty}}{\| h \|_{M,\infty}}.
}
Suppose that $E \subseteq \bbC$ is a compact set.  We now write $B(E)$ for the Banach space of functions that are continuous on $E$ and analytic in its interior with norm $\| f \|_{E} = \sup_{z \in E} | f(z) |$.

\thm{[\cite{TrefPlatteIllCond}]
\label{t:impossibility}
Let $E \subseteq \bbC$ be a compact set containing $[-1,1]$ in its interior and suppose that $\{ F_M \}^{\infty}_{M =1}$ is an approximation procedure based on equispaced grids of $M+1$ points such that for some $C,\rho > 1$ and $1/2 < \tau \leq 1$ we have
\bes{
\| f - F_M(f) \|_{\infty} \leq C \rho^{-M^{\tau}} \| f \|_{E},\qquad M=1,2,\ldots,
}
for all $f \in B(E)$.  Then the condition numbers \R{condition_number} satisfy
\bes{
\kappa(F_M) \geq \sigma^{M^{2\tau-1}},
}
for some $\sigma > 1$ and all large $M$.
}

Specializing to $\tau = 1$, this theorem states that exponential convergence at a geometric rate implies exponential ill-conditioning at a geometric rate.  Conversely, stability of any method is only possible when $\tau = 1/2$, which corresponds to root-exponential convergence in $M$.

\subsection{Coppersmith \& Rivlin's bound}
\label{ss:CRbound}
The proof of Theorem \ref{t:impossibility}, although it does not pertain to polynomials or polynomial approximation specifically, relies on a result of Coppersmith \& Rivlin on the maximal behaviour of polynomials bounded on an equispaced grid.  To state this result, we first introduce the following notation:
\be{
\label{BNMequispaced}
B(M,N) = \sup \left \{ \| p \|_{\infty} : p \in \bbP_{N}, \| p \|_{M,\infty} \leq 1 \right \}.
}
Note that in the special case $M = N$, this is just the Lebesgue constant
\bes{
B(N,N) = \Lambda(N) = \sup \left \{ \| F_N(f) \|_{\infty} : f \in C([-1,1]), \| f \|_{\infty} \leq 1 \right \},
}
where $F_N(f)$ denotes the polynomial interpolant of degree $N$ of a function $f$.

\thm{[\cite{copprivlinpolygrowth}]
\label{t:CoppRivOriginal}
Let $\{ x_m \}^{M}_{m=0} = \{ -1 + 2m/M \}^{M}_{m=0}$ be an equispaced grid of $M+1$ points in $[-1,1]$ and suppose that $1 \leq N \leq M$.  Then
there exist constants $b \geq a > 1$ such that
\bes{
a^{N^2/M} \leq B(M,N) \leq b^{N^2/M}.
}
}
Two implications of this result are as follows.  First, a polynomial of degree $N$ bounded on $M = \ord{N}$ equispaced points can grow at most exponentially large in between those points.  Second, one needs quadratically-many equispaced points, i.e.\ $M \asymp N^2$, in order to prohibit growth of an arbitrary polynomial of degree $N$ that is bounded on an equispaced grid.  We remark in passing that when $M = N$, so that $B(N,N) = \Lambda(N)$ is the Lebesgue constant, one also has the well-known estimate $\Lambda(N) \sim \frac{2^{N+1}}{\E N \log(N)}$ for large $N$ (see, for example, \cite[Chpt.\ 15]{TrefethenATAP}).

\rem{
\label{r:equispaced_sufficient}
Sufficiency of the scaling $M \asymp N^2$ is a much older result than Theorem \ref{t:CoppRivOriginal}, dating back to work Sch\"onhage \cite{schonhagepolyinterp} and Ehlich \& Zeller \cite{EhlichZeller1,EhlichZeller2} in the 1960s.  Ehlich also proved unboundedness of $B(M,N)$ if $M = o(N^2)$ as $N \rightarrow \infty$ \cite{EhlichPoly}. More recently, Rakhmanov \cite{RakhmanovPolyBds} has given a very precise analysis of not just $B(M,N)$ but also the pointwise quantity $B(M,N,x) = \sup \{ |p(x)| : p \in \bbP_N, \| p \|_{M,\infty} \leq 1 \}$ for $-1 \leq x \leq 1$.
}

\subsection{Discrete least squares}
A simple algorithm that attains the bounds implied by Theorem \ref{t:impossibility} is discrete least-squares fitting with polynomials:
\be{
\label{discreteLSequispaced}
F_{M,N}(f) = \underset{p \in \bbP_{N}}{\operatorname{argmin}} \sum^{M}_{m=0} | f(x_m) - p(x_m) |^2.
}
Here $N \leq M$ is a parameter which is chosen to ensure specific rates of convergence.  The following result determines the conditioning and convergence of this approximation (note that this result is valid for arbitrary sets of points, not just equispaced):

\prop{
\label{p:discLS}
Let $\{ x_m \}^{M}_{m=0}$ be a set of $M+1$ points in $[-1,1]$ and suppose that $1 \leq N \leq M$.  If $f \in C([-1,1])$ and $F_{M,N}(f)$ is as in \R{discreteLSequispaced} then the error
\bes{
\| f - F_{M,N}(f) \|_{\infty} \leq (1 + \kappa_{M,N} ) \inf_{p \in \bbP_{N}} \| f - p \|_{\infty},
}
where $\kappa_{M,N} = \kappa(F_{M,N})$ is the condition number of $F_{M,N}$.  Moreover,
\be{
\label{kappaBMN}
B(M,N) \leq \kappa_{M,N} \leq \sqrt{M+1} B(M,N),
}
where $B(M,N)$ is as in \R{BNMequispaced}.  Furthermore, if $M = N$, i.e.\ $F_{N} = F_{N,N}$ is the polynomial interpolant of degree $N$, then
\be{
\label{kappaLambdaN}
\kappa_{N,N} = B(N,N) = \Lambda(N),
}
is the Lebesgue constant.
}
Although this result is well known, we include a short proof for completeness:
\prf{
Since the points are distinct and $N \leq M$, the least-squares solution exists uniquely.  Notice that the mapping $F_{M,N}$ is linear and a projection onto $\bbP_{N}$.  Hence
\be{
\label{kappa_linearo}
\kappa_{M,N} = \sup_{\substack{f \in C([-1,1]) \\ \| f \|_{M,\infty} \neq 0}} \frac{\| F_{M,N}(f) \|_{\infty}}{\| f \|_{M,\infty}},
}
and consequently we have
\bes{
\| f - F_{M,N}(f) \|_{\infty} \leq \| f - p \|_{\infty} + \| F_{M,N}(f-p) \|_{\infty} \leq (1 + \kappa_{M,N} ) \| f - p \|_{\infty},\quad \forall p \in \bbP_N.
}
It remains to estimate the condition number.  Because $F_{M,N}(f)$ is a polynomial, it follows that
\be{
\label{kappa_linear}
\kappa_{M,N} \leq B(M,N) \sup_{\substack{f \in C([-1,1]) \\ \| f \|_{M,\infty} \neq 0}} \frac{\| F_{M,N}(f) \|_{M,\infty}}{\| f \|_{M,\infty}}.
}
Now observe that
\bes{
\| F_{M,N}(f) \|^2_{M,\infty} \leq \sum^{M}_{m=0} | F_{M,N}(f)(x_m) |^2 = \frac{M+1}{2} \| F_{M,N}(f) \|^2_{M,2}.
}
Since $F_{M,N}(f)$ is the solution of a discrete least-squares problem it is a projection with respect to the discrete semi-inner product $\ip{\cdot}{\cdot}_{M}$.  Hence $\| F_{M,N}(f) \|_{M,2} \leq \| f \|_{M,2} \leq \sqrt{2} \| f \|_{\infty}$.  Combining this with the previous estimate gives the upper bound $\kappa_{M,N} \leq \sqrt{M+1} B(M,N)$.  For the lower bound, we use \R{kappa_linearo} and the fact that $F_{M,N}$ is a projection to deduce that
\bes{
\kappa_{M,N} \geq \max_{\substack{p \in \bbP_{N} \\ \| p \|_{M,\infty} \neq 0}} \frac{\| p \|_{\infty}}{\| p\|_{M,\infty}} = B(M,N).
}
This completes the proof of \R{kappaBMN}.  For \R{kappaLambdaN} we merely use the definition of $\Lambda(N)$.
}

\subsection{Examples of nonequispaced points}\label{ss:examples}
To illustrate our main results proved later in the paper, we shall consider points $-1 = x_0 < x_1 < \ldots < x_M = 1$ that are equidistributed with respect to so-called \textit{modified Jacobi} weight functions.  These are defined as
\be{
\label{modJacobiWeight}
\mu(x) = g(x) (1-x)^{\alpha}(1+x)^{\beta},
}
where $\alpha,\beta > -1$ and $g \in L^\infty(-1,1)$ satisfies $c_1 \leq g(x) \leq c_2$ almost everywhere.  Throughout, we assume the normalization
\bes{
\int^{1}_{-1} \mu(x) \D x = 1,
}
in which case the points $\{x_m \}^{M}_{m=0}$ are defined implicitly by
\be{
\label{mu_equidistributed}
\frac{m}{M} = \int^{x_m}_{-1} \mu(x) \D x,\quad m=0,\ldots,M.
}
Ultraspherical weight functions are special cases of modified Jacobi weight functions.  They are defined as
\be{
\label{ultrasphericalWeigt}
\mu(x) = c (1-x^2)^{\alpha},\qquad c = \left ( \int^{1}_{-1} (1-x^2)^{\alpha} \D x \right )^{-1},
}
for $\alpha > -1$.  Within this subclass, we shall consider a number of specific examples:
\begin{itemize}
\item[(U)] ($\alpha = 0$) The uniform weight function $\mu(x) = \frac12$, corresponding to the equispaced points $x_m = -1 + 2 \frac{m}{M}$.
\item[(C1)] ($\alpha = -1/2$) The Chebyshev weight function of the first kind $\mu(x) = \frac{1}{\pi\sqrt{1-x^2}}$, corresponding to the Chebyshev points.  Note that these points are roughly equispaced near $x=0$ and quadratically spaced near the endpoints $x = \pm 1$.  That is,  $|x_{1} +1 |, |x_{M-1} - 1| = \ord{M^{-2}}$. 
\item[(C2)] ($\alpha = \frac12$) The Chebyshev weight function of the second kind $\mu(x) =\frac{2}{\pi} \sqrt{1-x^2}$.  Note that the corresponding points are roughly equispaced near $x = 0$, but are sparse near the endpoints.  In particular, $|x_{1} +1 |, |x_{M-1} - 1| = \ord{M^{-1/2}}$.  
\end{itemize}
Recall that for (U) one requires a quadratic scaling $M \asymp N^2$ to ensure stability.  Conversely, for (C1) any linear scaling of $M = c N$ with $c>1$ suffices (see Remark \ref{r:ChebyLS}).  Since the points (C2) are so poorly distributed near the endpoints, we expect, and it will turn out to be the case, that a more severe scaling than quadratic is required for stability in this case.  

We shall also consider two further examples:
\begin{itemize}
\item[(UC)] ($\alpha = -1/4$) The corresponding points cluster at the endpoints, although not quadratically.  Specifically, $|x_{1} +1 |, |x_{M-1} - 1| = \ord{M^{-4/3}}.$ 
\item[(OC)] ($\alpha = -3/4$) The corresponding points overcluster at the endpoints: $|x_{1} +1 |, |x_{M-1} - 1| = \ord{M^{-4}}$. 
\end{itemize}
We expect (UC) to require a superlinear scaling of $M$ with $N$ for stability, although not as severe as quadratic scaling as in the case of (U).  Conversely, in (OC) it transpires that linear scaling suffices, but unlike the case of (C1), the scaling factor $c$ (where $M/N = c$) must be sufficiently large. 

\begin{figure}
\begin{center}
\hspace{.5cm} C2  \hspace{3.2cm} UC \hspace{3.2cm} C1 \hspace{3.2cm} OC \\
 \includegraphics[width=16cm]{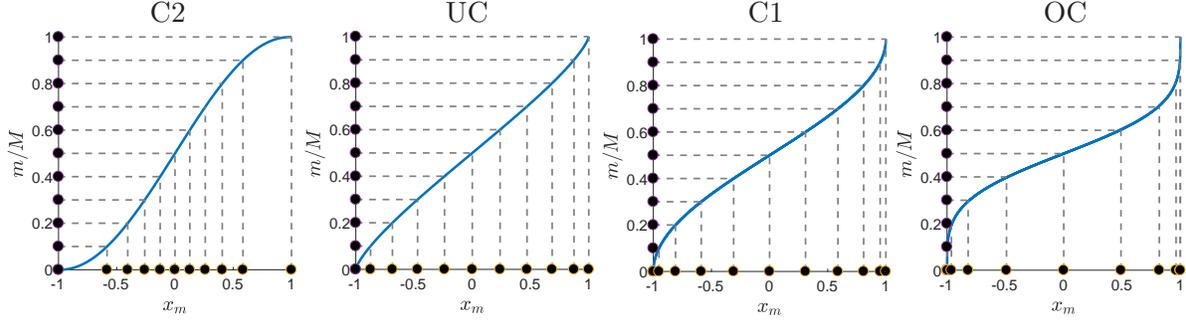}
  \end{center}
 \caption{Relationship between $m/M$ and $x_m$ given by (\ref{mu_equidistributed}) for the ultraspherical weight functions with $\alpha=1/2$ (C2), $\alpha=-1/4$ (UC), $\alpha=-1/2$ (C1), and $\alpha=-3/4$ (OC). Here $M=10$.}
\label{fig:points}
\end{figure}

The node clustering for the above distributions is illustrated in Figure~\ref{fig:points}.  This figure also shows the corresponding cumulative distribution functions $\int_{-1}^x \mu(\xi) d\xi$.

\section{Maximal behaviour of polynomials bounded on arbitrary grids}
We now seek to estimate the maximal behaviour of a polynomial of degree $N$ that is bounded at arbitrary nodes $-1 = x_0 < x_1 < \ldots < x_M = 1$.  As in \S \ref{ss:CRbound}, we define
\be{
\label{BMN}
B(M,N) = \sup \left \{ \| p \|_{\infty} : p \in \bbP_{N},\ | p(x_m) | \leq 1,\ m=0,\ldots,M \right \}.
}
Once again we note that $B(N,N) = \Lambda(N)$ is the Lebesgue constant of polynomial interpolation.

\subsection{Lower bound for $B(M,N)$}
\begin{figure}
\begin{center}
\begin{tabular}{c c}
$N=9$, $M=14$, $K_{-}=3$ & $N=15$, $M=30$, $K_{-}=4$  \\
\includegraphics[width=8cm]{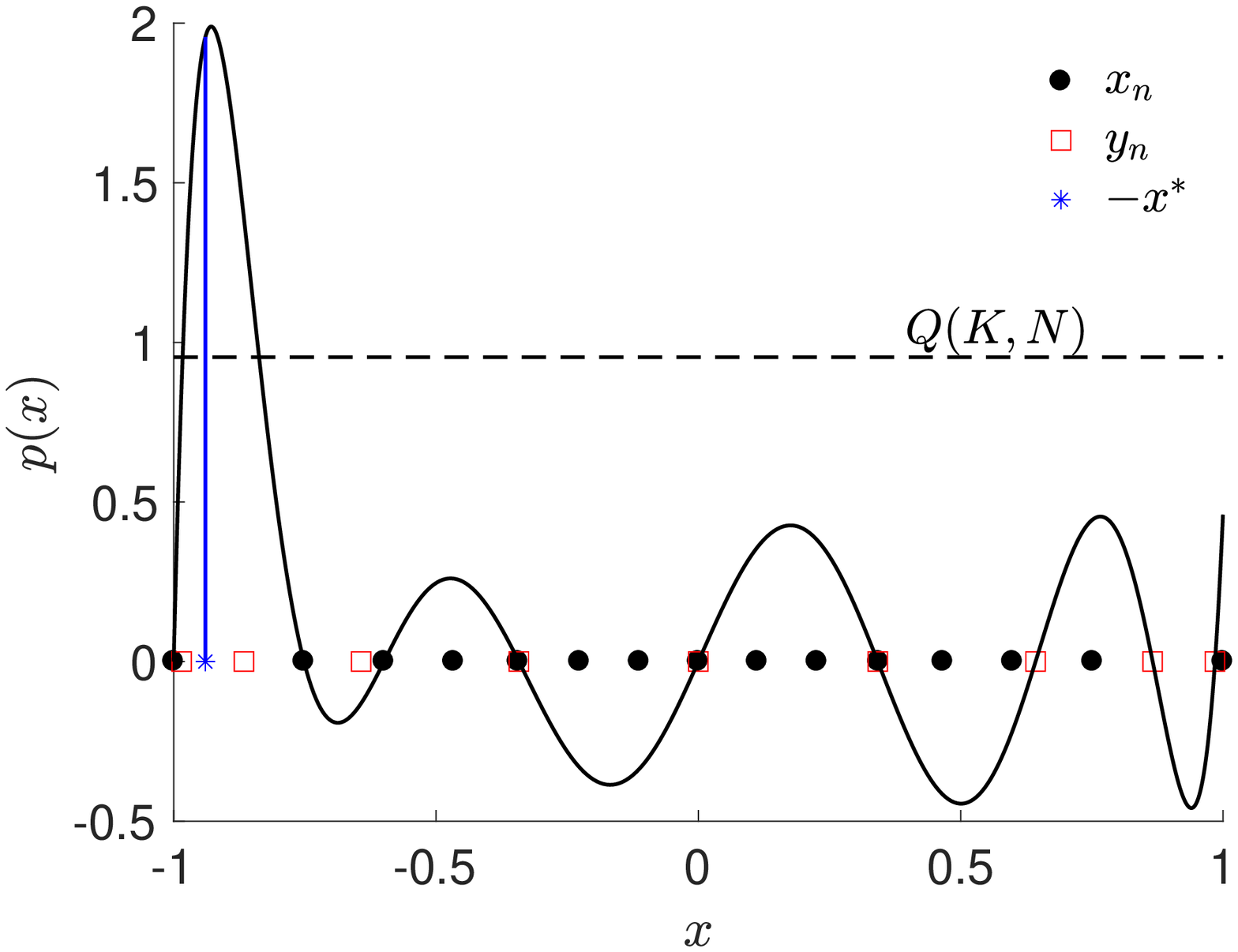} &\includegraphics[width=8cm]{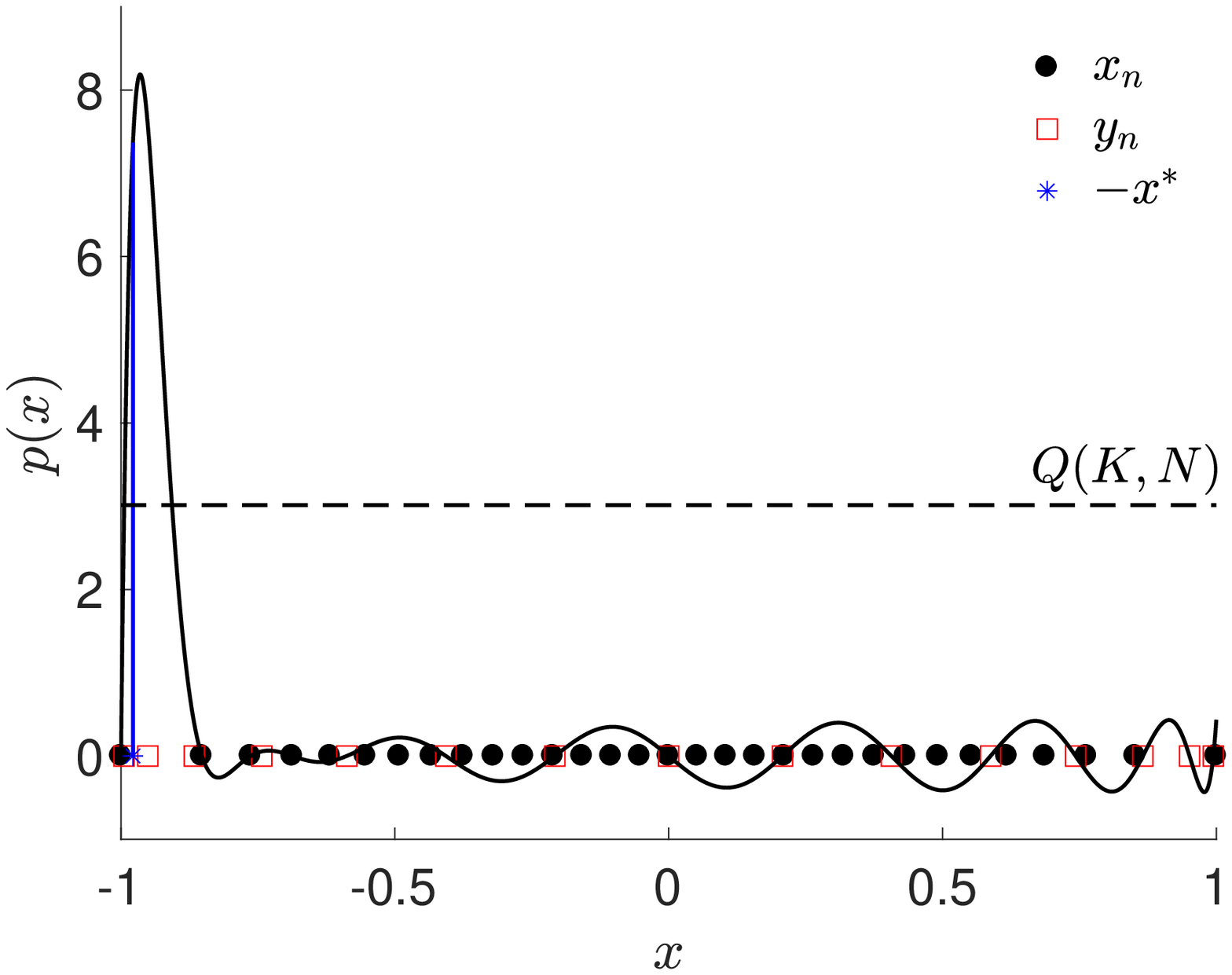} \\
\end{tabular}
\end{center}
\caption{Plots of the polynomial $p$ used in the proof of Theorem \ref{t:CoppRiv} for two sets of points. In both cases the points $x_m$ were generated with $\gamma=1/2$ -- case (C2) in \S \ref{ss:examples}. For reference, the lower bound $Q(K,N)$ is also included. }
\label{Fig_Theo_3_1}
\end{figure}

Our first main result of the paper is a generalization of the lower bound of Coppersmith \& Rivlin (Theorem \ref{t:CoppRivOriginal}) to arbitrary nodes.  Before stating this, we need several definitions.  First, given $M \geq N \geq 1$ and nodes $-1 = x_0 < x_1 < \ldots < x_M = 1$, we define 
\eas{
Q_{-}(K,N) &= \frac{\pi}{8}\left ( \frac{2N^2}{\pi^2} \right )^{K-1} \frac{1}{\Gamma(K+1/2)^2}  \prod^{K-1}_{n=1} (1+x_n),\qquad 2 \leq K \leq N,
\\
Q_{+}(K,N) &= \frac{\pi}{8}\left ( \frac{2N^2}{\pi^2} \right )^{K-1} \frac{1}{\Gamma(K+1/2)^2}  \prod^{K-1}_{n=1} (1-x_{M-n}),\qquad 2 \leq K \leq N,
\\
Q_{-}(1,N) &= Q_{+}(1,N) = 1.
}
Second, let $-1 < y_0 < \ldots < y_{N-1} < 1$ be the zeros of the $N^{\rth}$ Chebyshev polynomial $q(x) = \cos(N \arccos(x))$:
\be{
\label{Cheb_pts}
y_{n} = - \cos \left ( \frac{(2n+1) \pi}{2N} \right ),\qquad n=0,\ldots,N-1.
}
We now have the following:

\thm{
\label{t:CoppRiv}
Let $M \geq N \geq 1$, $-1 = x_0 < x_1 < \ldots < x_M = 1$ and suppose that there exist integers $2 \leq K_{-} \leq N$ and $2 \leq K_{+} \leq N$ such that
\be{
\label{K_def}
0 \geq x_{n} > y_{n},\qquad n=1,\ldots,K_{-}-1,
}
and
\be{
\label{K_def_2}
0 \leq x_{M-n} < y_{N-n},\qquad n=1,\ldots,K_{+}-1,
}
respectively, where the $y_n$ are as in \R{Cheb_pts}.  If either $K_-$ or $K_+$ fails to exist, set $K_- = 1$ or $K_+ = 1$.  Then the constant $B(M,N)$ defined in \R{BMN} satisfies
\be{
B(M,N) \geq \max \left \{ Q_{-}(K_{-},N), Q_{+}(K_+,N)\right \}.
}
}
\prf{
We first show that $B(M,N) \geq Q_{-}(K_{-},N)$.  If $K_{-} = 1$ there is nothing to prove, hence we now assume that $2 \leq K_{-} \leq N$.  Let $q(x) = \cos(N \arccos(x))$ be the $N^{\rth}$ Chebyshev polynomial and define $p \in \bbP_N$ by
\be{
\label{applepie}
p(x) = \frac12  q(x) \prod^{K_{-}-1}_{n=0} \frac{x-x_{n}}{x-y_n}.
}
Figure~\ref{Fig_Theo_3_1} illustrates the behaviour of $p$. We first claim that
\be{
\label{dinner}
|p(x_m) |\leq 1,\quad m=0,\ldots,M.
}
Clearly, for $m=0,\ldots,K_{-}-1$ we have $p(x_m) = 0$.  Suppose now that $K_{-} \leq m \leq M$.  Then, since $|q(x)| \leq 1$,
\bes{
| p(x_m) | \leq \frac12  \prod^{K_{-}-1}_{n=0} \left | \frac{x_m - x_{n}}{x_m - y_n} \right |.
}
By definition, we have $x_m > x_{n}$ for $n=0,\ldots,K_{-}-1$.  Also, by \R{K_def},
\bes{
x_{m} > x_{K_{-}-1} \geq y_{K_{-}-1} \geq y_{n},\quad n=0,\ldots,K_{-}-1,
}
and therefore
\bes{
|p(x_m) | \leq \frac12 \prod^{K_{-}-1}_{n=0} \frac{x_m - x_{n}}{x_m - y_n} .
}
For $n=1,\ldots,K_{-}-1$ \R{K_def} gives that $ \frac{x_m - x_{n}}{x_m - y_n} \leq 1$.  Also, since $x_m \geq y_{K_{-}-1}$ and $y_1 > - 1$ we have
\eas{
\frac{x_m - x_0}{x_m - y_0} \leq \frac{y_{K_{-}-1} +1}{y_{K_{-}-1}-y_0} 
 = 1 + \frac{1+y_0}{y_{K_{-}-1}-y_0}
 = 1 + \frac{\sin^2(\pi/(4N))}{\sin(K_{-} \pi/(2N)) \sin((K_{-}-1) \pi/(2N))}.
}
Recall that $2 t / \pi \leq \sin(t) \leq t$ for $0 \leq t \leq \pi/2$.  Hence
\bes{
\frac{x_m - x_0}{x_m - y_0} \leq 1 + \frac{\pi^2}{16 K_{-}(K_{-}-1)} \leq 2,
} 
and therefore
\bes{
|p(x_m) | \leq \frac12 \prod^{K_{-}-1}_{n=0} \frac{x_m-x_n}{x_m-y_n} \leq 1.
}
This completes the proof of the claim \R{dinner}.

We now wish to estimate $\| p \|_{\infty}$ from below.  Following Figure~\ref{Fig_Theo_3_1}, we choose the point $-x^* = -\cos ( \pi/N)$ midway between the endpoint $x=-1$ and the leftmost node $y_1$.  Since $|q(-x^*)|=1$ we derive from \R{applepie} that
\bes{
\| p \|_{\infty} \geq | p(-x^*) | = \frac12  \prod^{K_{-}-1}_{n=0} \left | \frac{x^* + x_n}{x^*+y_n} \right |.
}
Notice that $x^*+y_n > 0$ for $n=1,\ldots,K_{-}-1$ and therefore $x^*+x_n > 0$ for $n=1,\ldots,K_{-}-1$ by \R{K_def}.  Hence
\be{
\label{intermediate}
\| p \|_{\infty} \geq \frac12  \left | \frac{x^*+x_0}{x^*+y_0} \right | \prod^{K_{-}-1}_{n=1} \frac{x^*+x_n}{x^*+y_n} \geq \frac12  \left | \frac{x^*+x_0}{x^*+y_0} \right | \prod^{K_{-}-1}_{n=1} \frac{1+x_n}{1+y_n},
}
where in the second step we use \R{K_def} and the fact that $-y_n < x^* \leq 1$ and $x_n > y_n$.  Note that
\bes{
1+y_{n} = 2 \sin^2 \left ( \frac{(2n+1) \pi}{4N} \right )^2 \leq \frac{(2n-1)^2 \pi^2}{8 N^2},
}
and that
\bes{
\frac{|x^*+x_0|}{|x^*+y_0|} = \frac{1-x^*}{\cos(\pi/(2N)) - x^*} \geq 1.
}
Therefore we deduce that
\eas{
\| p \|_{\infty} &\geq \frac12   \left ( \prod^{K_{-}-1}_{n=1} \frac{8 N^2}{(2n-1)^2 \pi^2} \right ) \left ( \prod^{K_{-}-1}_{n=1} (1+x_n)  \right ) 
\\
&= \frac12  \left ( \frac{8 N^2}{\pi^2} \right )^{K_{-}-1} \frac{\pi}{4^K_{-} \Gamma(K_{-}+1/2)^2} \prod^{K_{-}-1}_{n=1}(1+x_n)
\\
& = Q_{-}(K_{-},N)
}
which gives $B(M,N) \geq Q_{-}(K_{-},N)$ as required.  In order to prove $B(M,N) \geq Q_{+}(K_{+},N)$ we repeat the same arguments, working from the right endpoint $x = +1$.
}

\begin{figure}
\begin{center}
(C2) \hspace{4.5cm} (U) \hspace{4.5cm} (UC)  \\
\hspace{-0.5cm} \includegraphics[width=5.5cm]{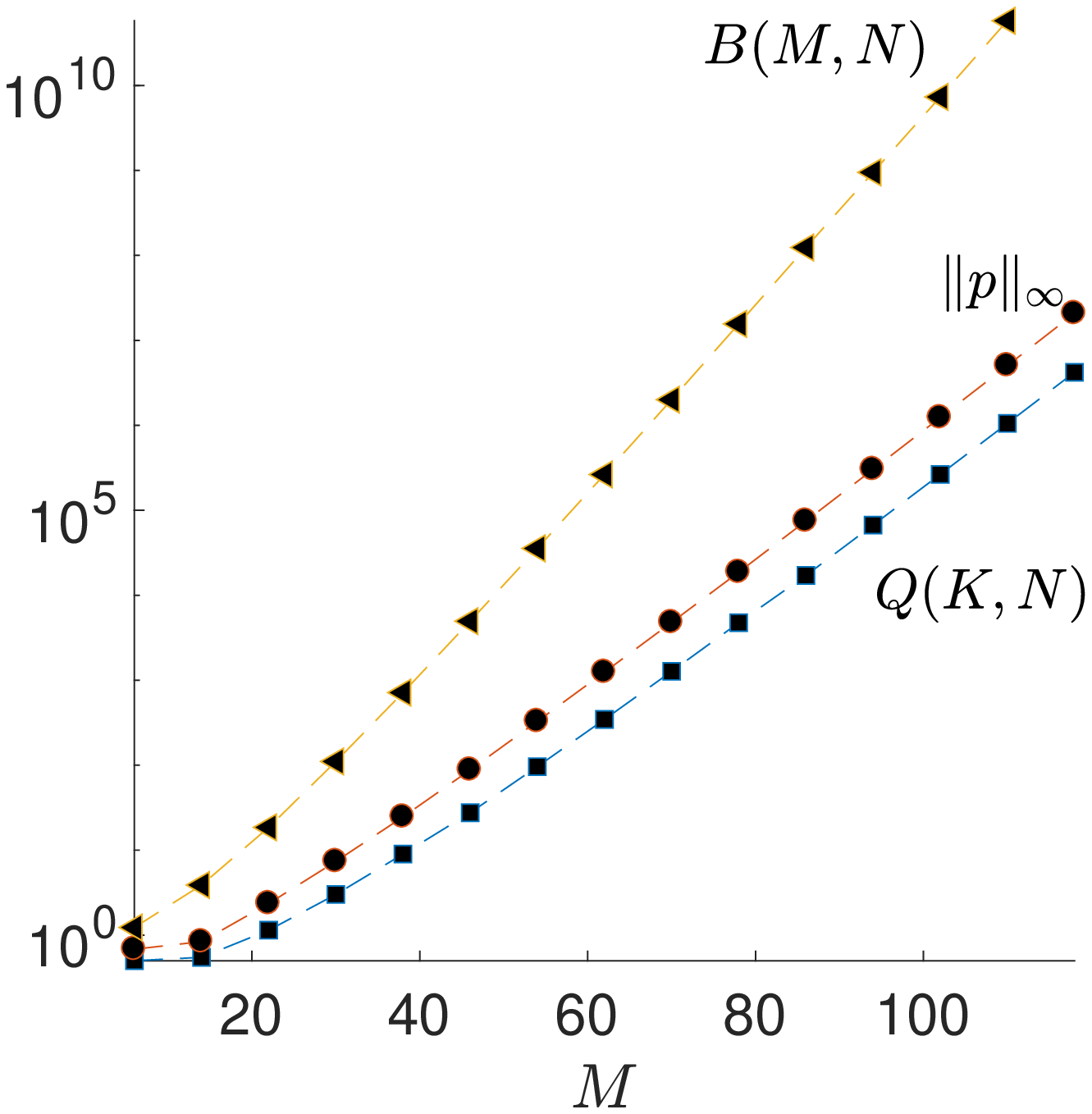} \hspace{-0.5cm} \includegraphics[width=5.5cm]{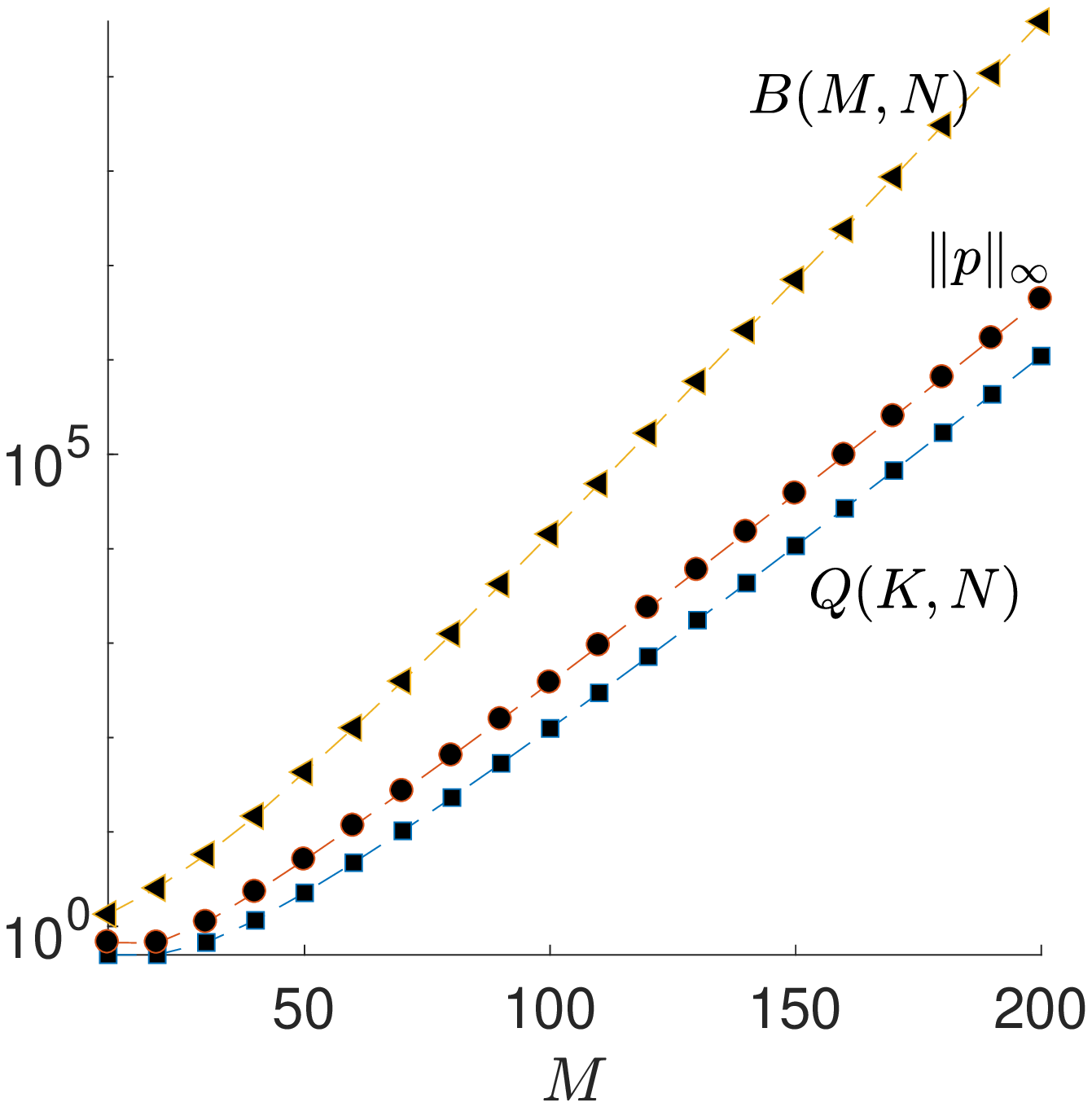}  \hspace{-0.5cm}\includegraphics[width=5.5cm]{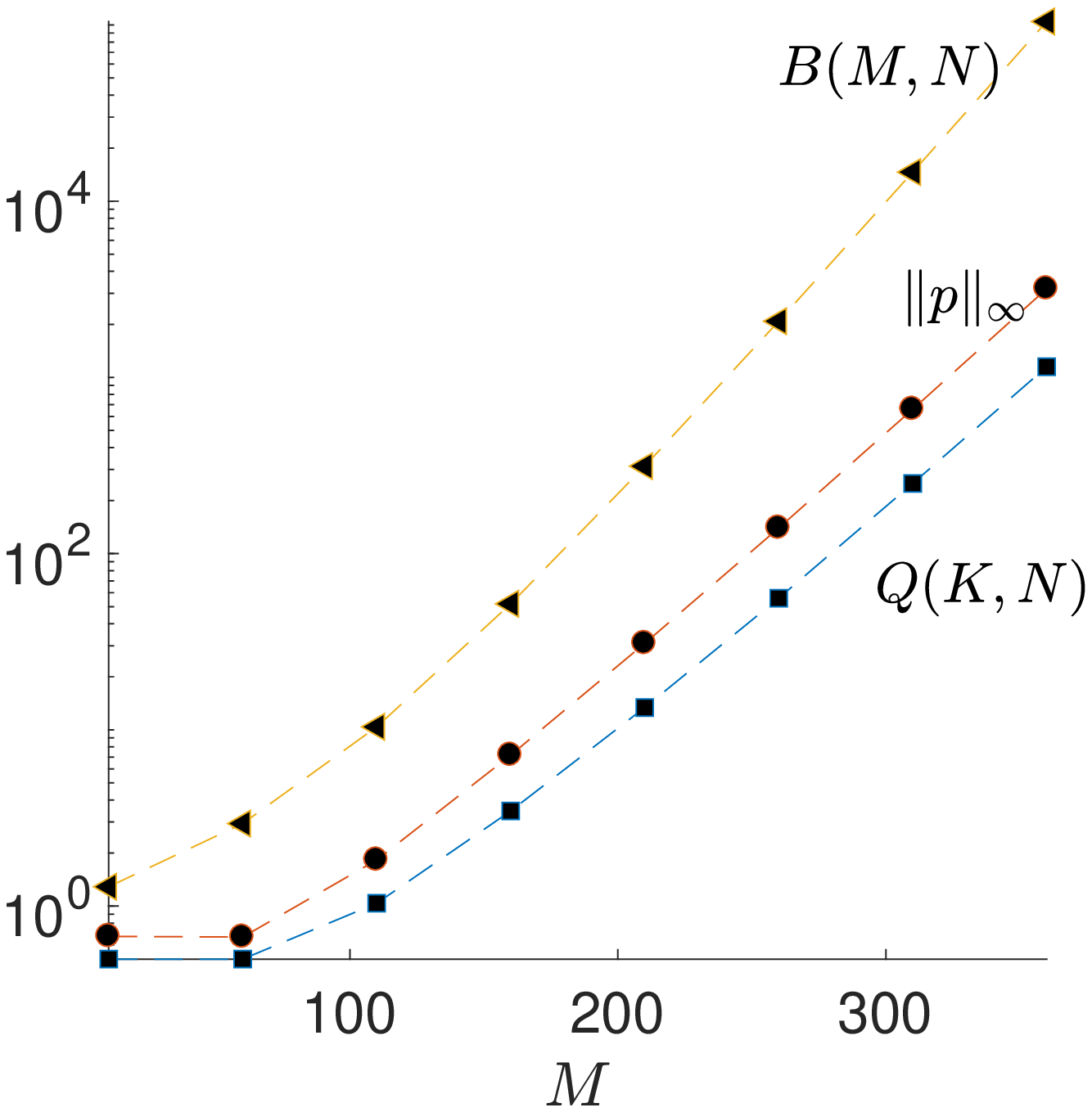}\\
\end{center}
\caption{Semi-log plot of the growth of $B(M,N)$, $\|p\|_\infty$, and $Q(K,N)$ for several values of $M$ and three node density functions described in \S\ref{ss:examples}: (C2), (U) and (UC). In all cases $N=M/2$. The computation of the quantity $B(M,N)$ is described in \S \ref{s:Remez}.}
\label{fig:Compare_BMN_QKN}
\end{figure}

Figure~\ref{fig:Compare_BMN_QKN} shows the growth of $B(M,N)$, $Q(K,N)$ and the norm of the polynomial used to prove Theorem \ref{t:CoppRiv}. In these examples, the nodes $x_m$ were generated using the density functions (C2), (U) and (UC). In all cases the polynomial degree was chosen as $N=M/2$. Notice that the exponential growth rate of $\|p\|_\infty$ is well estimated by $Q(K,N)$, while both quantities underestimate the rate of growth of $B(M,N)$.

\subsection{Lower bound for modified Jacobi weight functions}
Theorem \ref{t:CoppRiv} is valid for arbitrary sets of points $\{x_m \}^{M}_{m=0}$.  In order to derive rates of growth, we now consider points equidistributed according to modified Jacobi weight functions.  For this, we first recall the following bounds for the Gamma function (see, for example, \cite[Eqn.\ (6.1.38)]{AS}):
\be{
\label{Gamma_bounds}
\sqrt{2 \pi} z^{z+1/2} \E^{-z} \leq \Gamma(z+1) \leq \sqrt{2 \pi} \E^{1/12} z^{z+1/2} \E^{-z},\qquad z \geq 1.
}

\cor{
\label{c:modJacobi}
Let $\mu(x) = g(x) (1-x)^{\alpha} (1+x)^{\beta}$, where $\alpha , \beta > -1$ and $c_1 \leq g(x) \leq  c_2$ almost everywhere.  If $\gamma = \max \{ \alpha , \beta \} > -1/2$ then there exist constants $C > 0$ and $\sigma > 1$ depending on $\alpha$ and $\beta$ such that
\bes{
B(M,N) \geq C \sigma^{\nu},\qquad \nu = \left ( N^{2(\gamma+1)}/M \right )^{\frac{1}{2\gamma+1}},
}
for all $1 \leq N \leq M$.
}
\prf{
By definition, the points $\{x_m\}^{M}_{m=0}$ satisfy
\bes{
\frac{m}{M} = \int^{x_m}_{-1} g(x) (1-x)^{\alpha} (1+x)^{\beta} \D x.
}
Without loss of generality, suppose that $\beta \geq \alpha$.  Let $m$ be such that $x_m \leq 0$.  Then
\bes{
\frac{m}{M} = \int^{x_m}_{-1} g(x) (1-x)^{\alpha} (1+x)^{\beta} \D x \lesssim \int^{x_m}_{-1} (1+x)^{\beta} \D x \lesssim (1+x_m)^{1+\beta},
} 
and therefore
\be{
\label{xm_sing_meas}
1+x_m \gtrsim \left ( \frac{m}{M} \right )^{\frac{1}{1+\beta}},\qquad x_m \leq 0.
}
We now apply Theorem \ref{t:CoppRiv} and \R{Gamma_bounds} to get
\ea{
B(M,N) \geq Q(K,N) &\gtrsim \left ( \frac{2N^2}{\pi^2 M^{\frac{1}{1+\beta}}} \right )^{K-1} \frac{\Gamma(K)^{\frac{1}{1+\beta}} }{\Gamma(K+1/2)^2} \nn
\\
& \gtrsim  \left ( \frac{2N^2}{\pi^2 M^{\frac{1}{1+\beta}}} \right )^{K-1} \frac{\left((K/2)^{K-1} \E^{1-K}\right)^{\frac{1}{1+\beta}}}{K^{2K} \E^{2(1-K)}} \nn
\\
\label{BK_sing_meas}
& = K^{-2} \left ( \frac{2^{\frac{\beta}{1+\beta}} \E^{\frac{1+2 \beta}{1+\beta}}  N^2}{\pi^2 M^{\frac{1}{1+\beta}} K^{\frac{1+2\beta}{1+\beta}} } \right )^{K-1}
}
where $K$ is any value such that $x_n \geq y_n$ for $n=1,\ldots,K-1$.  We need to determine the range of $K$ for which this holds.  From \R{xm_sing_meas} we find that $x_n \geq y_n$ provided
\bes{
\left ( \frac{n}{M} \right )^{\frac{1}{1+\beta}} \gtrsim \frac{n^2}{N^2},
}
or equivalently
\bes{
n \leq c \left ( \frac{N^{2(\beta+1)}}{M} \right )^{\frac{1}{2\beta+1}} = c \nu,
}
where $c>0$ is some constant.  Hence \R{BK_sing_meas} holds for 
\be{
\label{Kregime}
K \in \{2,\ldots, 1 + \lfloor c \nu \rfloor \}.
}
We next pick a constant sufficiently small so that $0 < d \leq 2c/3$ and
\be{
\label{dcrit}
\frac{2^{\frac{\beta}{1+\beta}} \E^{\frac{1+2 \beta}{1+\beta}}  }{\pi^2  (2d)^{\frac{1+2\beta}{1+\beta}} } \geq 2.
}
Consider the case where
\bes{
\nu > 2 / d.
}
We set
\bes{
K = 1 + \left \lceil \frac{d}{2} \nu \right \rceil ,
}
and notice that
\bes{
K \leq 1 + 1+ \frac{d}{2} \nu < d \nu + \frac{d}{2} \nu = \frac{3 d}{2} \nu \leq c \nu.
}
 due to the assumptions on $d$.  Hence this value of $K$ satisfies \R{Kregime}.  We next apply \R{BK_sing_meas}, the bounds $K \geq d \nu / 2$ and $K \leq 2 d \nu$ and \R{dcrit} to deduce that
 \bes{
 B(M,N) \gtrsim \nu^{-2} \left ( 2^{d/2} \right )^{\nu} \geq \rho^{\nu},
 }
 for some $\rho > 1$.  This holds for all $\nu > 2/d$.  But since $d$ is a constant, we deduce that $B(M,N) \gtrsim \rho^{\nu}$ for all $\nu \geq 1$.  This completes the proof.
}

This result shows that if $M = \ord{N^{\tau}}$ for some $0 < \tau \leq 2(\gamma+1)$ then the maximal polynomial grows at least exponentially fast with rate
\bes{
r = \frac{2(\gamma+1) - \tau}{2\gamma+1}.
}
In particular, the scaling $M  \asymp N^{2(\gamma+1)}$, $N \rightarrow \infty$, is necessary for boundedness of the maximal polynomial.  In \S \ref{ss:nec_and_suffic} we will show that this rate is also sufficient.

It is informative to relate this result to several of the examples introduced in \S \ref{ss:examples}.  First, if $\gamma = 0$, i.e.\ case (U), we recover the lower bound of Theorem \ref{t:CoppRivOriginal}.  Conversely, for (C2) we have
\bes{
B(M,N) \geq C \sigma^{\sqrt{N^3/M}}.
}  
Thus, the cubic scaling $M \asymp N^3$ is necessary for stability.  Finally, for the points (UC) we have
\bes{
B(M,N) \geq C \sigma^{(M^{3/2}/N)^2},
}
which implies a necessary scaling of $M \asymp N^{3/2}$.  Note that Corollary \ref{c:modJacobi} says nothing about the case $-1 < \gamma \leq -1/2$.  We discuss the case $\gamma = -1/2$ further in Remark \ref{r:ChebyLS}.

The case of linear oversampling ($M \asymp N$) warrants closer inspection:
\cor{
\label{c:modJacobiLinOS}
Let $\mu(x)$, $C$ and $\sigma$ be as in Corollary \ref{c:modJacobi} and $c \geq 1$.  Then
\bes{
B(\lfloor cN \rfloor , N) \geq C \left ( \sigma^{c^{-\frac{1}{2\gamma+1}}} \right )^{N}.
}
In particular, the Lebesgue constants $\Lambda(N) = B(N,N)$ satisfy
\bes{
\Lambda(N) \geq C \sigma^{N}.
}
}

\begin{figure}
\begin{center}
\begin{tabular}{c c c}
(C2) & (UC) & (OC) \\
\includegraphics[width=5cm]{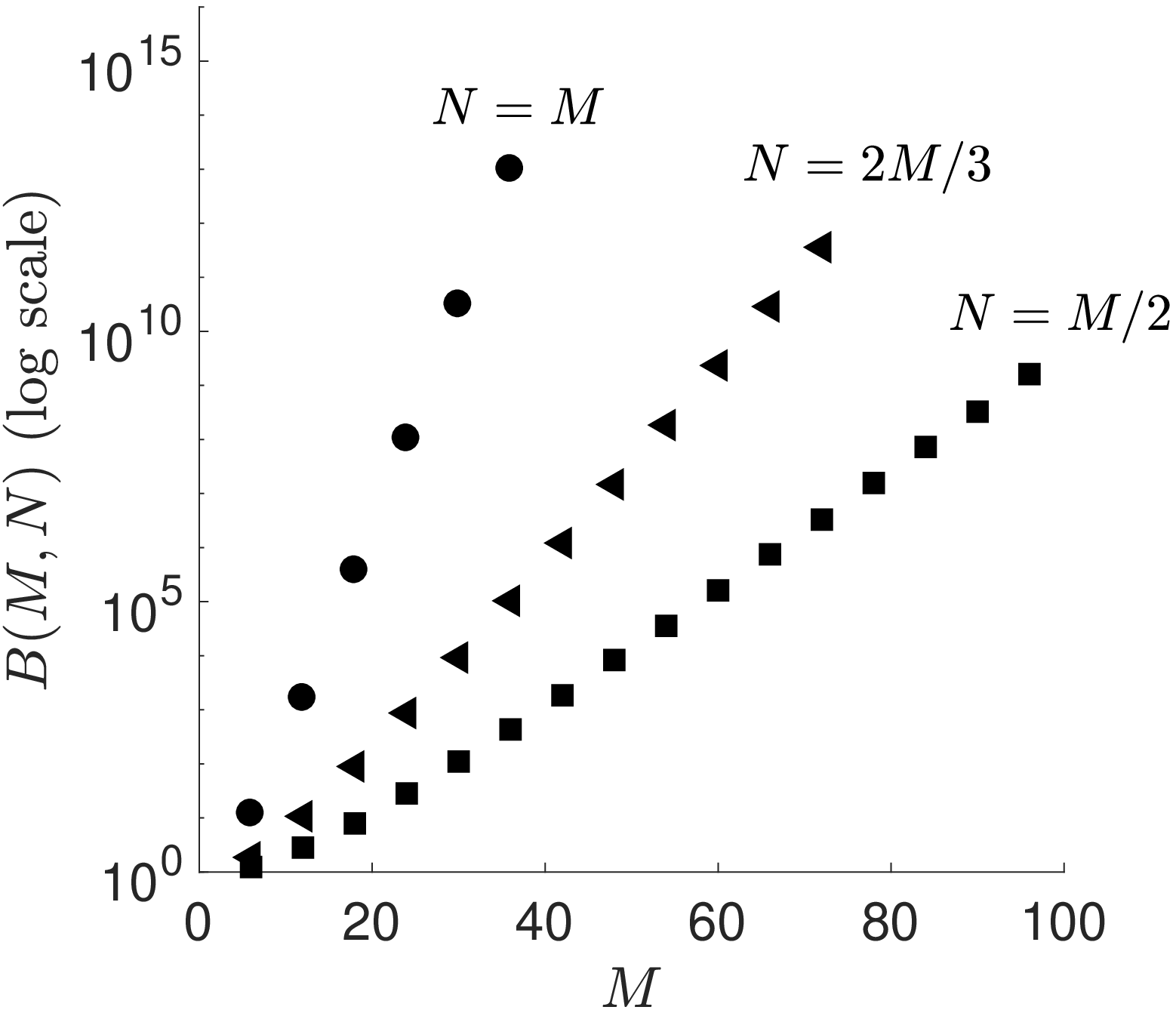} &
\includegraphics[width=5cm]{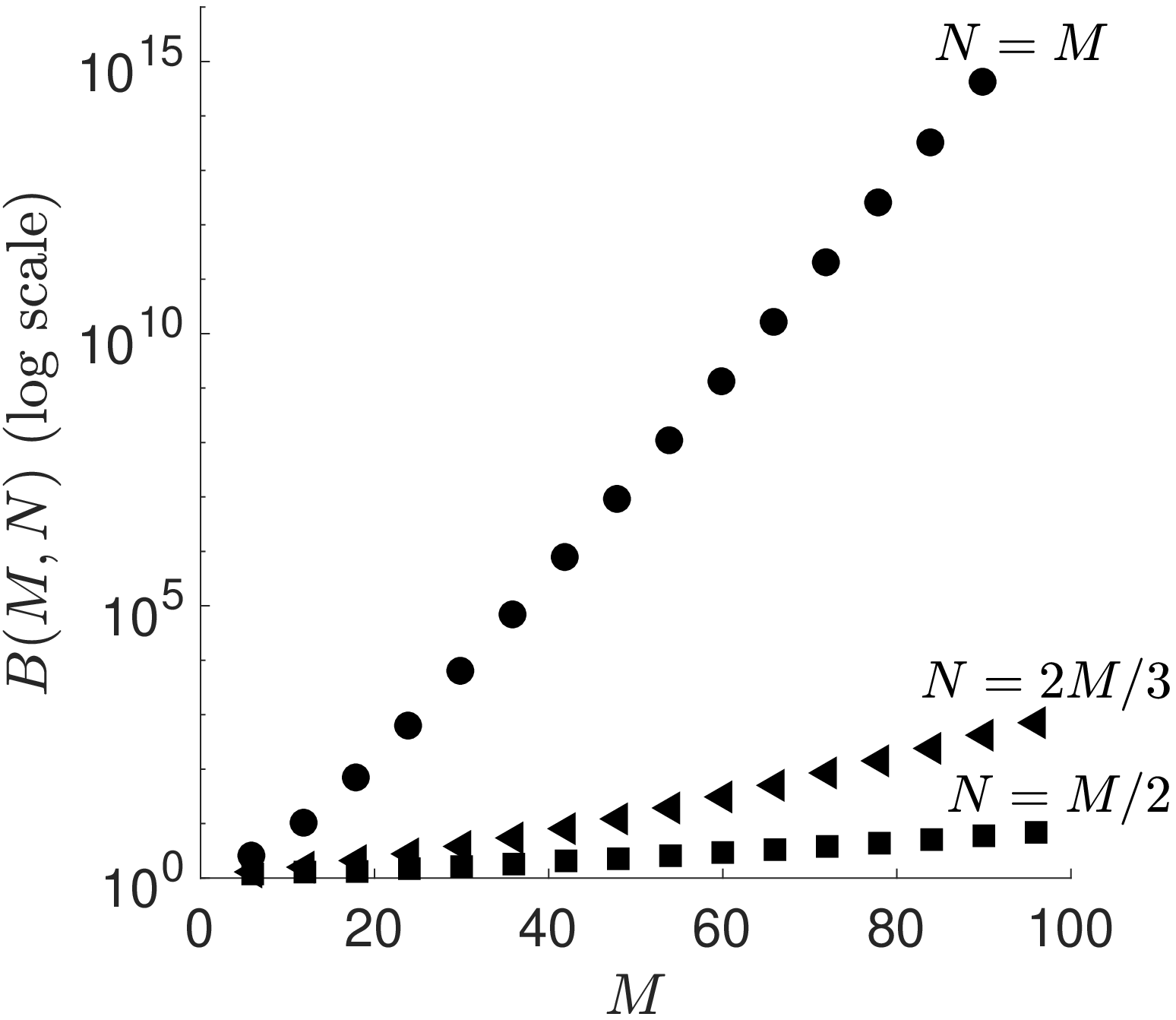}&\includegraphics[width=5cm]{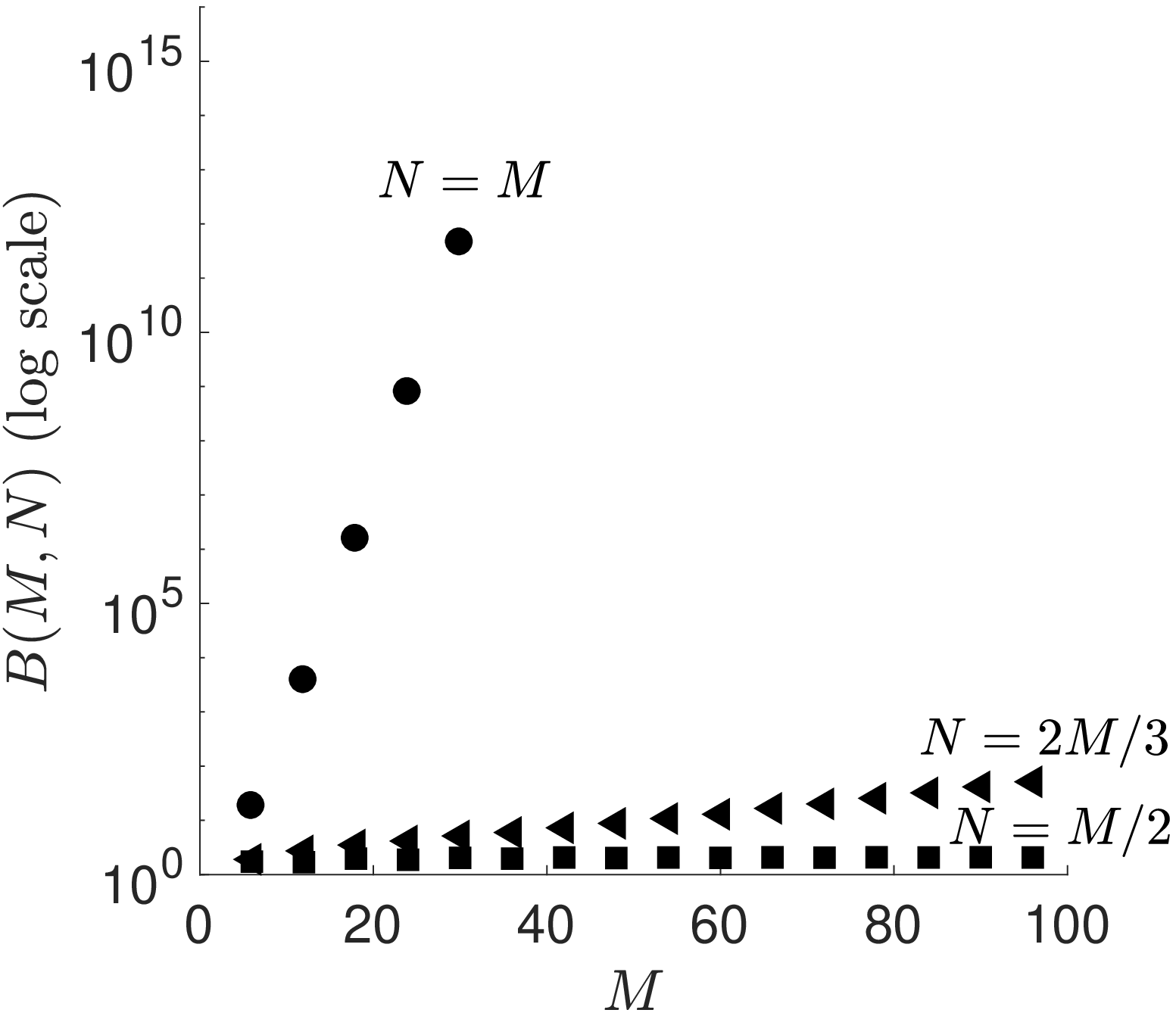}
\end{tabular}
\end{center}
\caption{The growth of $B(M,N)$ as a function of $M$ for the node densities (C2), (UC), and (OC). In each case $B(M,N)$ is plotted for $N = M$, $N = 2/3M$, and $N=1/2M$. }
\label{fig:growthB}
\end{figure}

In other words, whenever the points cluster more slowly than quadratically at one of the endpoints (recall that $\gamma = \max \{ \alpha ,\beta \} > -1/2$ in the above result), linear oversampling (including the case $M = N$, i.e.\ interpolation) necessarily leads to exponential growth of the maximal polynomial a geometric rate.  Figure~\ref{fig:growthB} illustrates this for the cases (C2) and (UC).  Interestingly, the case (OC), although not covered by this result, also exhibits exponential growth, whenever the oversampling factor is below a particular threshold.

\rem{
\label{r:Rakhmanov_upper}
In this paper we are primarily interested in lower bounds for $B(M,N)$, since this is all that is required for the various impossibility theorems.  However, in the case of ultraspherical weight functions an upper bound can be derived from results of Rakhmanov.  Specifically, let $\gamma = \alpha = \beta > -1/2$ and $x_1,\ldots,x_M$ be a set of $M$ points that are equispaced with respect to the ultraspherical weight function \R{ultrasphericalWeigt}.  Then \cite[Thm.\ 1(a)]{RakhmanovPolyBds} gives that
\bes{
\max_{-r \leq x \leq r} |p(x) | \leq C \max \{ |p(x_m) |, m=1,\ldots,M \},\qquad \forall p \in \bbP_N,
}
for some constant $C = C_r$ depending only on $r$, where $r^2 < 1 - (N/M)^{\tau}$ and $\tau = \frac{2}{2\gamma+1}$\footnote{Rakhmanov's result excludes the endpoints $x = \pm 1$ from this set, whereas in our results we include these points.  However as noted earlier, our main theorems would remain valid (with minor changes) if these points were excluded.}.  Remez' inequality (see, for example, \cite[Thm.\ 5.1.1]{borweinPoly}) now gives that
\bes{
\max_{|x| \leq 1} |p(x)| \leq C T_{N}(1/r) \max \{ |p(x_m) |, m=1,\ldots,M \},\qquad \forall p \in \bbP_N,
}
where $T_N$ is the $N^{\rth}$ Chebyshev polynomial.  Suppose that $N/M \leq 1/2$.  Then
\bes{
T_N(1/r) < \left ( 1/r + \sqrt{1/r^2-1} \right )^N < \left ( 1 + (N/M)^{\tau} + 2 (N/M)^{\tau/2} \right )^{N} \leq \exp(c (N/M)^{\tau/2} N),
}
for some $c>0$.  Using the definition of $\tau$, we obtain
\bes{
\| p \|_{\infty} \leq \tilde{C} \tilde{\sigma}^{\nu} \max \{ |p(x_m) |, m=1,\ldots,M \},\qquad \forall p \in \bbP_N,
}
for $\tilde{C} > 0$ and $\tilde{\sigma} > 1$, where
\bes{
\nu = \left ( \frac{N^{2(\gamma+1)}}{M} \right )^{\frac{1}{2\gamma+1}}.
}
The exponent $\nu$ is exactly the same as in the lower bound for $B(M,N)$ (Corollary \ref{c:modJacobi}).  In other words, the two-sided bounds of Coppersmith \& Rivlin (see Theorem \ref{t:CoppRivOriginal}) for equispaced nodes extend to nodes equidistributed according to ultraspherical weight functions.
}

\rem{
\label{r:ChebyLS}
Corollaries \ref{c:modJacobi} and \ref{c:modJacobiLinOS} do not apply to the nodes (C1).  It is, however, well-known (see, for example, \cite[Thm.\ 15.2]{TrefethenATAP}) that
\bes{
\Lambda(N) \sim \frac{2}{\pi} \log(N),\quad N \rightarrow \infty,
}
in this case.  Furthermore, a classical result of Ehlich \& Zeller \cite{EhlichZeller1} gives that
\bes{
B^*(M,N) \leq \frac{1}{\cos \left ( \frac{\pi N}{2M} \right )},
}
where
\bes{
B^*(M,N) = \sup\left \{ \| p \|_{\infty} : p \in \bbP_N , |p(y_m)| \leq 1,\ m=1,\ldots,M \right \},
}
and $y_m = \cos \left ( \frac{2m-1}{2M} \right )$, $m=1,\ldots,M$.\footnote{Similar to the previous footnote, $B^*(M,N)$ excludes the endpoints $x = \pm 1$.}  In particular, this result implies boundedess of $B^*(M,N)$ in the case of linear oversampling, i.e.\ $M \geq (1+\epsilon) N$ for any $\epsilon > 0$.  
}

\section{An extended impossibility theorem}

For $\theta > 1$, let $E_{\theta} \subseteq \bbC$ be the Bernstein ellipse with parameter $\theta$.  That is, the ellipse with foci at $\pm 1$ and semiminor and semimajor axis lengths summing to $\theta$.  Given a domain $E \subseteq \bbC$, we let $D(E)$ be the set of functions that are analytic on $E$.  The following lemma -- whose proof follows the same ideas to that of Theorem \ref{t:impossibility} -- shows how the condition number of an exponentially-convergent method for approximating analytic functions can be bounded in terms of the quantity $B(M,N)$ for suitable $N$.

\lem{[Abstract impossibility lemma]
\label{l:absimpossibility}
Given points $-1 = x_0 < x_1 < \ldots < x_M = 1$, let $F_M : C([-1,1]) \rightarrow C([-1,1])$ be  an approximation procedure such that $F_M(f)$ depends only on the values $\{ f(x_m) \}^{M}_{m=0}$ for any $f \in C([-1,1])$.  Suppose that
\be{
\label{FMdecay}
\| f - F_M(f) \|_{\infty} \leq C \rho^{-M^{\tau}} \| f \|_{E},
}
for all $f \in D(E)$, where $E \subseteq \bbC$ be a compact set containing $[-1,1]$ in its interior, and $C>0$, $\rho > 1$ and $\tau > 0$ are constants that are independent of $f$ and $M$.  If $N \in \bbN$ satisfies
\bes{
N < \frac{M^{\tau} \log(\rho) - \log(2C)}{\log(\theta)},
}
where $\theta > 1$ is such that the $E \subseteq E_{\theta}$, then the condition number $\kappa(F_M)$ defined in \R{condition_number} satisfies
\bes{
\kappa(F_M) \geq \frac12 B(M,N),
}
for $B(M,N)$ is as in \R{BMN}.
}
\prf{
Let $p \in \bbP_{N}$.  Since $p$ is entire \R{FMdecay} gives
\bes{
\| F_M(p) \|_{\infty} \geq \| p \|_{\infty} - C \rho^{-M^{\tau}} \| p\|_{E}.
}
Also, due to a classical result of Bernstein \cite[Sec.\ 9]{BernsteinPolyBound} (see also \cite[Lem.\ 1]{TrefPlatteIllCond}), one has $\| p \|_{E} \leq \| p \|_{E_{\theta}} \leq \theta^{N} \| p \|_{\infty}$.  Hence
\bes{
\| F_M(p) \|_{\infty} \geq \left ( 1 - C \theta^N \rho^{-M^{\tau}} \right ) \| p \|_{\infty} \geq \frac12 \| p \|_{\infty}.
}
It now follows that
\bes{
\kappa(F_M) \geq \max_{\substack{p \in \bbP_{N-1} \\ \| p \|_{M,\infty} \neq 0}} \frac{\| F_M(p) \|_{\infty}}{\| p \|_{M,\infty}} \geq \frac12 \max_{\substack{p \in \bbP_{N-1} \\ \| p \|_{M,\infty} \neq 0}}  \frac{\| p \|_{\infty}}{\| p \|_{M,\infty}}= \frac12 B(M,N),
}
as required.
}

In order to demonstrate this result in a concrete setting, we specialize to the case of points equidistributed with respect to modified Jacobi weight functions.  This leads to the following theorem, which is the second main result of the paper:

\thm{[Impossibility theorem for modified Jacobi weight functions]
\label{t:extendedimpossibility}
For $M \in \bbN$, let $\{ x_m \}^{M}_{m=0}$ be equidistributed according to a modified Jacobi weight function \R{modJacobiWeight} with parameters $\alpha , \beta > -1$ (see \R{mu_equidistributed}). Let $F_M : C([-1,1]) \rightarrow C([-1,1])$ be a family of approximation procedures such that $F_M(f)$ depends only on the values $\{ f(x_m) \}^{M}_{m=0}$ for any $f \in C([-1,1])$ and $M \in \bbN$.  Suppose that
\be{
\label{FMdecay2}
\| f - F_M(f) \|_{\infty} \leq C \rho^{-M^{\tau}} \| f \|_{E},
}
for all $f \in D(E)$, where $E \subset \bbC$ is a compact set containing $[-1,1]$ in its interior, and $C>0$, $\rho > 1$ and $\tau > 0$ are independent of $f$ and $M$.  If
\bes{
\gamma = \max \{ \alpha , \beta \} > -1/2,
}
and
\bes{
\tau > \frac{1}{2(\gamma+1)},
}
then the condition numbers $\kappa(F_M)$ satisfy
\bes{
\kappa(F_M) \geq \sigma^{M^{\nu}},
}
for some $\sigma > 1$ and all large $M$, where
\bes{
\nu = \frac{2(\gamma+1) \tau - 1}{2 \gamma+1} .
}
}
\prf{
We apply Lemma \ref{l:absimpossibility} and Corollary \ref{c:modJacobi}.
}

This result is best summarized by the statements in following corollary:

\cor{
\label{c:imposs_conclusions}
Consider the setup of Theorem \ref{t:extendedimpossibility}.  If $\gamma = \max \{ \alpha,\beta \} > -1/2$ then the following holds:
\begin{itemize}
\item[(i)] If $F_M(f)$ converges exponentially fast with geometric rate for all $f \in B(E)$ (i.e.\ $\tau = 1$ in \R{FMdecay2}) then the condition numbers $\kappa(F_M) \geq \sigma^{M}$ grow exponentially fast with geometric rate.
\item[(ii)] The best possible rate of exponential convergence of a stable method $F_M$ is subgeometric with index $\frac{1}{2(\gamma+1)}$.
\end{itemize}
}
\prf{
We set $\tau = 1$ (part (i)) or $\tau = \frac{1}{2(\gamma+1)}$ (part (ii)) in Theorem \ref{t:extendedimpossibility}.
}

Note that by letting $\gamma = 0$ (i.e.\ equispaced points) we recover the original impossibility theorem (Theorem \ref{t:impossibility}).  It is of interest to note that geometric convergence necessarily implies geometrically large condition numbers, regardless of the endpoint behaviour of the sampling points whenever $\gamma = \max \{ \alpha,\beta \} > -1/2$.  Conversely, this result says nothing about points that cluster quadratically or faster at $x = \pm 1$, which corresponds to the case $-1 < \gamma \leq -1/2$.  Indeed we shall prove in the next section that geometric convergence is possible in this setting with a stable approximation.

\section{Optimality of the approximation rate and discrete least squares}
Theorem \ref{t:extendedimpossibility} gives a necessary relation between the rate of exponential convergence and the rate of exponential ill-conditioning.  For example, as asserted in Corollary \ref{c:imposs_conclusions}, stable approximation necessarily implies subgeometric convergence with index $\frac{1}{2(\gamma+1)}$.  We now show that there exists an algorithm that achieves these rates: namely, polynomial least-squares fitting.  In particular, if the polynomial degree is chosen as
\bes{
N \asymp M^{\nu},\qquad \nu = \frac{1}{2(\gamma+1)},
}
one obtains a stable approximation which converges exponentially with rate $\frac{1}{2(\gamma+1)}$.

\subsection{A sufficient condition for boundedness of the maximal polynomial}\label{ss:nec_and_suffic}

We commence with a sufficient condition for boundedness of the quantity $B(M,N)$:

\lem{
\label{l:B_zeta}
Let $-1 = x_0 < \ldots < x_M = 1$ be arbitrary and suppose that $N \zeta < 1$, where
\be{
\label{zeta_def}
\zeta = \max_{m=0,\ldots,M-1} \int^{x_{m+1}}_{x_m} \frac{1}{\sqrt{1-x^2}} \D x.
}
If $B(M,N)$ is as in \R{BMN}, then
\bes{
B(M,N) \leq \frac{1}{1-N \zeta}.
}
}
\prf{
Let $p \in \bbP_{N}$ with $| p(x_m) | \leq 1$, $m=0,\ldots,M$ and suppose that $-1 \leq x \leq 1$ with $x_{m} \leq x \leq x_{m+1}$ for some $m=0,\ldots,M-1$.  Then
\bes{
|p(x)| \leq | p(x_m) | + \int^{x}_{x_m} | p'(z) | \D z \leq 1 + \sup_{-1 \leq z \leq 1} | \sqrt{1-z^2} p'(z) | \int^{x_{m+1}}_{x_m} \frac{1}{\sqrt{1-z^2}} \D z.
}
Bernstein's inequality states that $| \sqrt{1-z^2} p'(z) | \leq N \| p \|_{\infty}$ for any $-1 \leq z \leq 1$ \cite[Thm.\ 5.1.7]{borweinPoly}.  Hence
\bes{
| p(x) | \leq 1 + N \zeta \| p \|_{\infty}.
}
Since $-1 \leq x \leq 1$ was arbitrary the result now follows.
}

Note that this lemma makes no assumption on the points $\{ x_m \}^{M}_{m=0}$.  The following result estimates the constant $\zeta$ for points arising from modified Jacobi weight functions:

\lem{
\label{l:modJacobi_suffic}
Let $\mu$ be a modified Jacobi weight function \R{modJacobiWeight} with parameters $\alpha , \beta > -1$. If $\zeta$ is as in \R{zeta_def}, then 
\bes{
\zeta \leq C M^{-\frac{1}{2(1+\gamma)}},
}
where $\gamma = \max \{ \alpha , \beta , -1/2 \}$ and $C>0$ is a constant depending on $\alpha$ and $\beta$ only.
}
\prf{
We consider the following four cases:
\bull{
\item[(i)] $-1 < \alpha , \beta \leq -1/2$,
\item[(ii)] $-1 < \alpha \leq -1/2$, $\beta > -1/2$,
\item[(iii)] $\alpha > -1/2$, $-1 < \beta \leq -1/2$,
\item[(iv)] $\alpha , \beta > - 1/2$.
}
Case (i): Suppose first that $-1 < \alpha , \beta \leq -1/2$.  Then
\bes{
\int^{x_{m+1}}_{x_m} \frac{1}{\sqrt{1-x^2}} \D x \lesssim \int^{x_{m+1}}_{x_m} \mu(x) \D x = \frac{1}{M}.
}
Hence $\zeta \lesssim 1/M$ in this case.  

\pbk
Case (ii): Now suppose that $-1 < \alpha \leq -1/2$ and $\beta > -1/2$.  Then
\bes{
\int^{x_{m+1}}_{x_m} \frac{1}{\sqrt{1-x^2}} \D x \lesssim \int^{x_{m+1}}_{x_m} g(x) (1-x)^{\alpha} (1+x)^{-1/2} \D x.
}
Recall from \R{xm_sing_meas} that
\bes{
1+x_m \gtrsim \left ( \frac{m}{M} \right )^{\frac{1}{1+\beta}},\qquad m=0,\ldots,M,
}
and therefore $1+x_{1} \gtrsim M^{-\frac{1}{1+\beta}}$.  For $m =1,\ldots,M-1$ we  have
\eas{
\int^{x_{m+1}}_{x_m} g(x) (1-x)^{\alpha} (1+x)^{-1/2} \D x &\leq (1+x_m)^{-1/2 - \beta} \int^{x_{m+1}}_{x_m} g(x) (1-x)^{\alpha} (1+x)^{\beta} \D x 
\\
& =  (1+x_m)^{-1/2 - \beta} \frac{1}{M}
\\
& \leq (1+x_1)^{-1/2-\beta} \frac{1}{M}
\\
& \lesssim M^{-\frac{1}{2(1+\beta)}}.
}
Now let $m = 0$.  For this, we first notice that $x_{1} \leq 0$ whenever $M \gtrsim 1$.  Therefore, we have
\bes{
\int^{x_1}_{x_0} g(x)(1-x)^{\alpha}(1+x)^{-1/2} \D x \lesssim \sqrt{1+x_1} \lesssim M^{-\frac{1}{2(1+\beta)}}.
}
Hence we deduce that $\zeta \lesssim M^{-\frac{1}{2(1+\beta)}}$ for this case as well.

\pbk
Case (iii): This is identical to the previous case and thus omitted.

\pbk
Case (iv): For $m=1,\ldots,M-2$ we have
\eas{
\int^{x_{m+1}}_{x_m} \frac{1}{\sqrt{1-x^2}} \D x &\leq (1-x_{m+1})^{-1/2-\alpha} (1+x_m)^{-1/2-\beta} \int^{x_{m+1}}_{x_m} g(x) (1-x)^{\alpha}(1+x)^{\beta} \D x 
\\
& \leq (1-x_{M-1})^{-1/2-\alpha} (1+x_1)^{-1/2-\beta} \frac{1}{M}
\\
& \lesssim M^{1-\frac{1}{2(1+\alpha)} -\frac{1}{2(1+\beta)}}
\\
& \leq M^{-\frac{1}{2(1+\gamma)}},
}
where in the final step we use the fact that $\alpha,\beta > -1/2$.
For $m=0$, recalling that $x_1 \leq 0$ for all large $M$, we have
\bes{
\int^{x_1}_{x_0} \frac{1}{\sqrt{1-x^2}} \D x \lesssim \sqrt{1+x_1} \lesssim M^{-\frac{1}{2(1+\beta)}},
}
and similarly for $m=M-1$, noting that $x_{M-1} \geq 0$ for all large $M$ gives
\bes{
\int^{x_M}_{x_{M-1}} \frac{1}{\sqrt{1-x^2}} \D x \lesssim \sqrt{1-x_M} \lesssim M^{-\frac{1}{2(1+\alpha)}}.
}
We therefore deduce that $\zeta \lesssim M^{-\frac{1}{2(1+\gamma)}}$ as required.
}

This lemma immediately leads to the follow result:

\prop{[Necessary and sufficient condition for boundedness of $B(M,N)$]
\label{p:nec_suffic_B}
For $N,M \in \bbN$, let $\{ x_m \}^{M}_{m=0}$ be equidistributed according to a modified Jacobi weight function \R{modJacobiWeight} with parameters $\alpha , \beta > -1$.  Then $B(M,N) \lesssim 1$ if and only if
\bes{
M \asymp N^{2(\gamma+1)}.
}
}
\prf{
We combine Corollary \ref{c:modJacobi} and Lemma \ref{l:modJacobi_suffic}.
}

\begin{figure}
\centerline{\includegraphics[width=.6\textwidth]{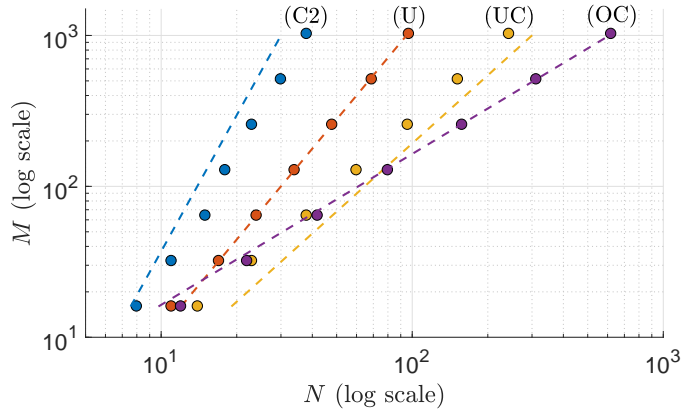}}
\caption{Loglog plot of $N$ and $M$, where $M$ is the smallest integer such that  $B(M,N)\le 10$. Circle markers show the computed values and dashed lines represent the theoretical estimate in Proposition~\ref{p:nec_suffic_B}. The dashed lines corresponding to (U), (C2) and (UC) are given by  $M = (N/3)^{2(\gamma+1)}$, while the dashed line corresponding to (OC) is $M = 1.65N$. }
\label{fig:prop5_3}
\end{figure}

\begin{figure}
\hspace{3cm} (U) \hspace{3.7cm} (UC) \hspace{3.8cm} (OC) \\
\includegraphics[width=1.1\textwidth]{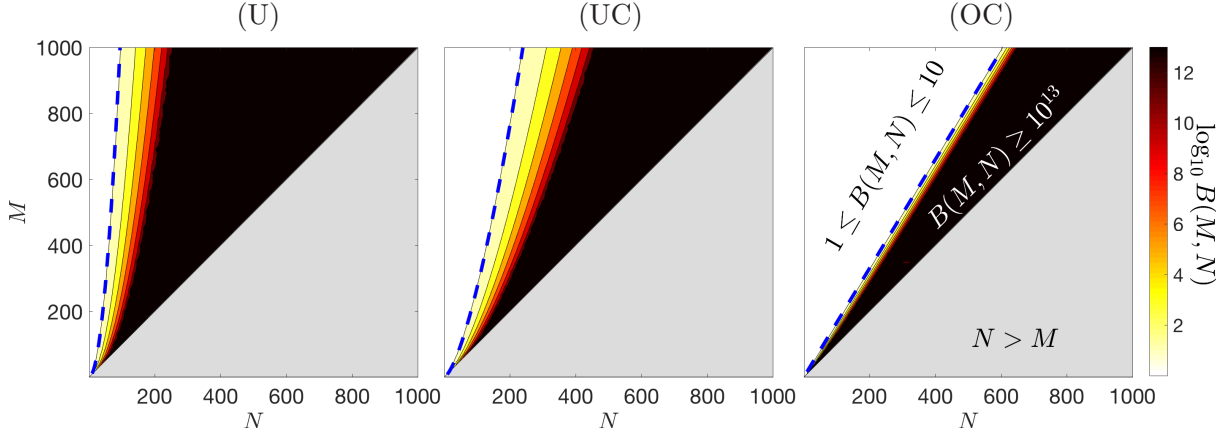}
\caption{Contour plot of $\log_{10} B(M,N)$ for node densities (U), (UC), and (OC). Black regions correspond to $B(M,N)>10^{13}$ and white regions to $B(M,N)<10$.  Dashed lines are given by $M=(N/3)^2$ (U), $M = (N/2.4)^{3/2}$ (UC), and $M = 1.65N$ (OC).}
\label{fig:prop5_3B}
\end{figure}

The proposition is illustrated in Figure~\ref{fig:prop5_3} for cases (C2), (U), (UC) and (OC). It plots the smallest values of $M$ such that $B(M,N)\le10$ for given values of $N$. Notice that the relationship between the computed values of $M$ and $N$ is in good agreement with the asymptotic relation $M \asymp N^{2(\gamma+1)}$. The constants used to define the dashed lines in this figure were chosen by trial and error. Similar agreement is shown in Figure~\ref{fig:prop5_3B}, where the contour levels of $B(M,N)$ for cases (U), (UC), and (OC) are presented. Notice that for the (OC) case $B(M,N)$ remains bounded with $M=CN$, but its values very quickly increase from 10 to more than $10^{13}$ when $C$ is decreased below $1.65$.


\rem{
For equispaced points, the sufficiency of the rate $M \asymp N^2$ is a classical result (see Remark \ref{r:equispaced_sufficient}).  More recently, similar sufficient conditions have appeared when the sampling points are drawn randomly and independently according the measure $\mu(x) \D x$.  For example, \cite{DavenportEtAlLeastSquares} proves that $M \asymp (N \log(N))^2$ uniformly-distributed points are sufficient for $L^2/\ell^2$-stability (note that we consider $L^{\infty}/\ell^\infty$-stability in this paper), whereas only $M \asymp N \log(N)$ points are required when drawn from the Chebyshev distribution.  In the multivariate setting, similar results have been proved in \cite{MiglioratiThesis,MiglioratiEtAlFoCM} for quasi-uniform measures.  Up to the log factors and the different norms used, these are the same as the rate prescribed in Proposition \ref{p:nec_suffic_B}, which is both sufficient and necessary.
}

\subsection{Application to polynomial least squares}\label{ss:LSapplic}
We now apply Proposition \ref{p:nec_suffic_B} to show that polynomial least squares achieves the optimal approximation rate of a stable approximation (up to a small algebraic factor in $M$) specified by the generalized impossibility theorem (Theorem \ref{t:extendedimpossibility}):

\prop{
\label{p:least_squares_optimal}
For $M \in \bbN$, let $\{ x_m \}^{M}_{m=0}$ be equidistributed according to a modified Jacobi weight function \R{modJacobiWeight} with parameters $\alpha , \beta > -1$.  For $1 \leq N \leq M$, let $F_{N,M}$ be the discrete least-squares approximation defined by \R{discreteLSequispaced}.  If
\bes{
M \asymp N^{2(\gamma+1)},
}
then
\bes{
1 \leq \kappa(F_{M,N}) \leq C \sqrt{M},
}
and
\bes{
\| f - F_{N,M}(f) \|_{\infty} \leq \frac{C \sqrt{M}}{\theta - 1} \theta^{-M^{\frac{1}{2(\gamma+1)}} } \| f \|_{E_{\theta}},
}
for all $f \in D(E_{\theta})$ and any $\theta > 1$, where $C >0$ is a constant.
}
\prf{
Proposition \ref{p:discLS} gives $\kappa(F_{M,N}) \leq \sqrt{2} \sqrt{M} B(M,N)$ and
\bes{
\| f - F_{M,N}(f) \|_{\infty} \leq 2 \sqrt{2} \sqrt{M} B(M,N) \inf_{p \in \bbP_{N}} \| f - p \|_{\infty}.
}
The result follows immediately from Proposition \ref{p:nec_suffic_B} and the well-known error bound $\inf_{p \in \bbP_N} \| f - p \|_{\infty} \leq \frac{2}{\theta - 1} \| f \|_{E_{\theta}} \theta^{-N}$ (see, for example, \cite[Thm.\ 8.2]{TrefethenATAP}).
}

\subsection{The mock-Chebyshev heuristic}
A well-known heuristic is that stable approximation from a set of $M+1$ points $-1 = x_0 < x_1 < \ldots < x_M = 1$ is only possible once there exists a subset of $N+1$ of those points that mimic a Chebyshev grid (see, for example, \cite{boyd2009divergence}).  We now confirm this heuristic.  Let
\be{
\label{ChebNodes}
z_n = - \cos(n \pi/N),\qquad n=0,\ldots,N,
}
be a Chebyshev grid (note that the $z_n$ are equidistributed according to the Chebyshev weight function $\mu(x) = \frac{1}{\pi\sqrt{1-x^2}}$).  Then we have:

\prop{
\label{p:interlace}
Let $-1 = x_0 < \ldots < x_M = 1$ be arbitrary and suppose that
\bes{
N \zeta < \pi,
}
where $\zeta$ is as in \R{zeta_def}.  Then there exists a subset $\{ x^*_n \}^{N}_{n=1}$ of $\{ x_m \}^{M}_{m=0}$ such that 
\be{
\label{interlace}
-1 = z_0 < x^*_1 < z_1 < x^*_2 < \ldots < z_{N-1} < x^*_N < z_N = 1,
}
where the $z_n$ are as in \R{ChebNodes}.  In particular, if $\{ x_m \}^{M}_{m=0}$ are equidistributed according to a modified Jacobi weight function with parameters $\alpha,\beta > -1$ then there exists such a subset $\{ x^*_n \}^{N}_{n=1}$ whenever
\bes{
M \asymp N^{2(\gamma+1)},\qquad \gamma = \max \{ \alpha, \beta , -1/2 \}.
}
}
\prf{
Let $\theta_m = \cos^{-1}(-x_m) \in [0,\pi]$ so that
\be{
\label{}
\theta_{m+1} - \theta_{m} = \int^{x_{m+1}}_{x_m} \frac{1}{\sqrt{1-x^2}} \D x \leq \zeta < \pi/N.
}
We now construct a subset $\{ \phi_n \}^{N}_{n=1} \subseteq \{ \theta_m \}^{M}_{m=0}$ such that
\bes{
0 < \phi_1 < \frac{\pi}{N} < \phi_2 < \frac{2 \pi}{N} < \ldots < \frac{(N-1)\pi}{N} < \phi_N < \pi.
}
First, let $m_1$ be the largest $m$ such that $\theta_m < \pi/N$.  Set $\phi_1 = \theta_{m_1}$.  Next, observe that $\theta_{m_1+1} < \theta_{m_1} + \pi/N < 2 \pi/N$.  Hence there exists at least one of the $\theta_m$'s in the interval $(\pi/N,2 \pi/N)$.  Let $m_2$ be the largest $m$ such that $\theta_m < 2 \pi/N$ and set $\phi_2 = \theta_{m_2}$.  We now continue in the same way to construct a sequence $\{ \phi_n \}^{N}_{n=1}$ with the required property.  Since the function $\cos^{-1}(-\theta)$ is increasing on $[0,\pi]$ it follows that the sequence $\{ x^*_n \}^{N}_{n=1}$ with $x^*_n = \cos^{-1}(-\phi_n)$ satisfies \R{interlace}.
}

Recalling Lemma \ref{l:B_zeta}, we note that the same sufficient condition (up to a small change in the right-hand side) for boundedness of the maximal polynomial (for arbitrary points) also guarantees an interlacing property of a subset of $N$ of those points with the Chebyshev nodes.  In particular, if $M,N \rightarrow \infty$ with $N \zeta \leq c$, where $c < \pi$, the nodes $\{ x^*_n \}^{N}_{n=1}$ equidistribute according to the Chebyshev weight function $\mu(x) = \frac{1}{\pi\sqrt{1-x^2}}$.  Moreover, for modified Jacobi weight functions the sampling rate that guarantees the existence of this `mock-Chebyshev' grid, i.e.\ $M \asymp N^{2(\gamma+1)}$, is identical to that which was found to be both necessary and sufficient for stable approximation via discrete least-squares (recall Proposition \ref{p:least_squares_optimal}).

\section{Computation of $B(M,N)$ and the maximal polynomial}\label{s:Remez}
Let $-1 = x_0 < x_1 < \ldots < x_M = 1$ be a set of $M+1$ points.  We now introduce an algorithm for the computation of 
\bes{
B(M,N) = \max \left \{ \nm{p}_{\infty} : p \in \bbP_N , \ | p(x_m) | \leq 1,\ m=0,\ldots,M \right \},
}
and the maximizing polynomial $p \in \bbP_N$.  In fact, in order to compute this polynomial we will first compute the pointwise quantity
\be{
\label{BMNx}
B(M,N,x) = \max \left \{ | p(x) | : | p(x_m) | \leq 1,\ m=0,\ldots,M,\ p \in \bbP_{N} \right \},\quad -1 \leq x \leq 1.
}
As we prove below, $B(M,N,x)$ is a piecewise polynomial with knots at the points $\{ x_m \}^{M}_{m=0}$.  Hence the maximal polynomial for $B(M,N)$ can be obtained by computing $B(M,N,x)$ in each subinterval and identifying the interval and corresponding polynomial in which the maximum is attained.

Our algorithm for computing \R{BMNx} is a variant of the classical Remez algorithm for computing best uniform approximations (see, for example, \cite{PachonTrefethenRemez,powell}).  It is based on a result of Sch\"onhage \cite{schonhagepolyinterp}.

\subsection{Derivation}
We first require some notation.  Given a set $Y = \{ y_n \}^{N}_{n=0}$ of $N+1$ points, let 
\bes{
\ell_{Y,n}(x) = \prod^{N}_{\substack{m = 0 \\ m \neq n}} \frac{x-y_m}{y_n - y_m},\qquad n=0,\ldots,N,
}
be the Lagrange polynomials, and 
\bes{
L_{Y}(x) = \max \left \{ | p(x) | : p \in \bbP_N , | p(y_n) | \leq 1, n=0,\ldots,N \right \} = \sum^{N}_{n=0} | \ell_{Y,n}(x) |,
}
be the  Lebesgue function of $Y$.  Note the second equality is a straightforward exercise.  We now require the following lemma:

\lem{
\label{l:LagrangeEquiosc}
Let $y_0 < y_1 < \ldots < y_N$.  If $x \in [y_n , y_{n+1}]$ then
\bes{
L_{Y}(x) = p_{Y,n}(x),
}
where $p_{Y,n} \in \bbP_{N}$ is the unique polynomial such that
\be{
\label{LagrangeInterpPoly}
p_{Y,n}(y_k) = \left \{ \begin{array}{ll} (-1)^{n-k} & k = 0,\ldots,n \\ (-1)^{n+1-k} & k = n+1,\ldots,N \end{array} \right ..
}
}
\prf{
Notice that, for $x \in [y_n,y_{n+1}]$, we have
\bes{
\sgn \left ( \ell_{Y,k}(x) \right ) = \left \{ \begin{array}{ll} (-1)^{n-k} & k=0,\ldots,n \\ (-1)^{n+1-k} & k=n+1,\ldots,N \end{array} \right ..
}
Hence
\bes{
L_Y(x) = \sum^{N}_{k=0} | \ell_{Y,k}(x) | = \sum^{N}_{k=0} \sgn \left ( \ell_{Y,k}(x) \right ) \ell_{Y,k}(x) = \sum^{N}_{k=0} p_{Y,n}(y_k) \ell_{Y,k}(x) = p_{Y,n}(x),
}
as required.
}

\begin{figure}
\centerline{\includegraphics[width=.55\textwidth]{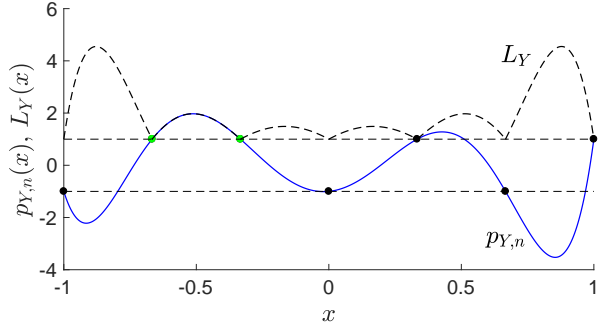}}
\caption{The Lebsgue function $L_{Y}$ for 7 equispaced points and its polynomial representation $p_{Y,n}$ on the interval $[y_1,y_2]$, as defined in Lemma~\ref{l:LagrangeEquiosc}.} 
\label{fig:LebLemma}
\end{figure}

This lemma is illustrated in Figure~\ref{fig:LebLemma}, where $L_Y$ and its polynomial representation $p_{Y,n}$ on the interval $[y_1,y_2]$ are plotted. 

\lem{
\label{l:technical_Rolle}
Let $y_0 < y_1 < \ldots < y_N$ and for $n=0,\ldots,N-1$ consider the polynomial $p_{Y,n}(x)$ defined in Lemma \ref{l:LagrangeEquiosc}.  Then
\bes{
p_{Y,n}(x) < 1,\quad x \in (y_{n-1},y_n) \cup (y_{n+1},y_{n+2}).
}
}
\prf{
Write $p = p_{Y,n}$ and $I_k = [y_k , y_{k+1}]$ for $k=0,\ldots,N-1$.  Since $p(y) > 1$ for $y \in I_n \backslash \{ y_n,y_{n+1}\}$ and $p(y_n) = p(y_{n+1}) = 1$ there is a point in $I_n$ where $p'$ vanishes.  Additionally, $p'(y_{n+1}) < 0$.  Suppose now that $p(y) \geq 1$ for some $y \in I_{n+1}$.  Then the derivative $p'$ will have at least two zeros in $I_{n+1}$.  The polynomial $p$ has at least $N-1$ zeros on $(y_0 , y_N)$, and since $p(y_n) = p(y_{n+1}) = 1$ it has no zeros on $I_n$.  Hence it must have exactly one zero on each subinterval $I_k$ for $k \neq n$ and one further zero outside $(y_0,y_N)$.  This implies there are at least $N-1$ zeros of $p$ outside $I_{n} \cup I_{n+1}$, and therefore $p'$ has at least $N-3$ zeros outside $I_{n} \cup I_{n+1}$.  Adding the one zero in $I_n$ and the two zeros in $I_{n+1}$ implies that $p'$ has at least $N$ zeros.  Since $p \in \bbP_N$ this is impossible.
}

We now produce our main result that will lead to the Remez-type algorithm.  The following result is due to Sch\"onhage \cite{schonhagepolyinterp}.  Since \cite{schonhagepolyinterp} is written in German and the relevant result (``Satz 3'') is stated for equispaced points only, we reproduce the proof below:

\lem{
\label{l:schonhage}
Let $-1 = x_0 < x_1 < \ldots < x_M = 1$ and $B(M,N,x)$ be as in \R{BMNx}.  If $x_m < x < x_{m+1}$ for some $m=0,\ldots,M-1$ then
\bes{
B(M,N,x) = \min_{Y} L_{Y}(x),
}
where the minimum is taken over all sets $Y \subseteq \{ x_m \}^{M}_{m=0}$ of size $|Y| =N+1$ with $x_m , x_{m+1} \in Y$.
}
\prf{
Consider a set $Y \subseteq \{ x_m \}^{M}_{m=0}$ of size $|Y| =N+1$ with $x_m , x_{m+1} \in Y$.  Then, by definition, $B(M,N,x) \leq L_Y(x)$.  Since there are only finitely-many such $Y$, there is a $Y^*$ with
\bes{
L_{Y^*}(x) = \min_{Y} L_Y(x).
}
Let $Y^* = \{ y_k \}^{N}_{k=0}$ and $p = p_{Y^*,n}$ be the polynomial defined in Lemma \ref{l:LagrangeEquiosc}, where $n$ is such that $y_n = x_m$.  We now claim that
\be{
\label{schonhage_claim}
| p(x_j) | \leq 1,\quad j=0,\ldots,M.
}
We shall prove this claim in a moment, but let us first note that this implies the lemma.  Indeed, assuming \R{schonhage_claim} holds we have
\bes{
p(x) \leq B(M,N,x) \leq L_{Y^*}(x) = p(x).
}
Hence $B(M,N,x) = p(x) = L_{Y^*}(x)$ as required.

To prove the claim we argue by contradiction.  Suppose that \R{schonhage_claim} does not hold and let $j$ be such that $| p(x_j) | > 1$.  Note that $x_j \notin Y^*$.  There are now three cases:

\pbk
Case 1: Suppose $x_j$ lies between two adjacent points of $Y^*$, i.e.\
\be{
\label{xj_between_yk}
y_k < x_j < y_{k+1},
}
for some $k = 0 ,\ldots,N-1$.  Since $\sgn(p(y_{k+1})) = - \sgn(p(y_k))$ there are two subcases:
\bes{
\mbox{(a):}\ \sgn(p(x_j)) = \sgn(p(y_k)),\qquad \mbox{(b):}\ \sgn(p(x_j)) = \sgn(p(y_{k+1})).
}
Suppose that subcase (a) occurs.  Exchange $y_k$ with $x_j$ and define the new set
\bes{
\widehat{Y} = (Y^* \backslash \{y_k \}) \cup \{ x_j \} = \{ \hat{y}_i \}^{N}_{i=0}.
}
We now claim that $x_m , x_{m+1} \in \widehat{Y}$; in other words, $j \neq m,m+1$.  First, notice that $j \neq m$ since \R{xj_between_yk} cannot hold when $k = n$.  Second, $j = m+1$ cannot hold either.  Indeed, if $j = m+1$ then $\sgn(p(x_j)) = \sgn(p(x_{m+1})) = 1$ and hence $p(x_j) > 1$.  But by Lemma \ref{l:technical_Rolle} we have $p(x) < 1$ for $y_{n+1} < x < y_{n+2}$, which is a contradiction.

Let $\ell_{i}$ and $\hat{\ell}_{i}$ be the Lagrange polynomials for $Y^*$ and $\widehat{Y}$ respectively.  Then, for $x_m < x < x_{m+1}$, we have
\be{
\sgn \left ( \hat{\ell}_i(x) \right ) = \sgn \left ( \ell_{i}(x) \right )= \sgn \left( p(y_i) \right ) = \sgn \left ( p(\hat{y}_i) \right ) .
}
Hence, expanding $p$ in the Lagrange polynomials $\hat{\ell}_i$, we obtain
\bes{
L_{Y^*}(x) = p(x) = \sum^{N}_{i=0} p(\hat{y}_i) \hat{\ell}_i(x) = \sum^{N}_{i=0} | p(\hat{y}_i) | | \hat{\ell}_i(x) | > \sum^{N}_{i=0} | \hat{\ell}_i(x) | = L_{\widehat{Y}}(x),
}
which contradicts the minimality of $Y^*$.  Subcase (b) is treated in a similar manner.

\pbk
Case 2: Suppose that $x_j > y_N$.  In this case we have the two subcases
\bes{
\mbox{(a):}\ \sgn(p(x_j)) = \sgn(p(y_N)),\qquad \mbox{(b):}\ \sgn(p(x_j)) =- \sgn(p(y_{N})).
}
In subcase (a) we construct $\hat{Y}$ by replacing $y_N$ with $x_j$ and, similarly to Case 1, arrive at a contradiction.

Now consider subcase (b).  We first note that $x_m \neq y_0$.  Indeed, if this were the case then $p$ would have $N+1$ zeros -- namely, $N-1$ zeros between $y_1$ and $y_N$, one zero between $y_N$ and $x_j$ and one zero to the left of $y_0$ -- which is a contradiction.  Hence we can exchange $y_0$ with $x_j$ to construct a new set of points
\bes{
\widehat{Y} = \left ( Y^* \backslash \{ y_0 \}  \right ) \cup \{ x_j \} = \{ y_{i} \}^{N+1}_{i=1},
}
where $y_{N+1} = x_j$.  The Lagrange polynomials on $\widehat{Y}$ satisfy
\bes{
\sgn \left ( \hat{\ell}_i(x) \right )  = \sgn \left ( \ell_i(x) \right ) = \sgn \left ( p(y_i) \right ),\quad i=1,\ldots,N,
}
and
\bes{
\sgn \left ( \hat{\ell}_{N+1}(x) \right ) =  - \sgn \left ( \ell_{N}(x) \right ) = - \sgn \left ( p(y_N) \right ) =\sgn \left ( p(y_{N+1}) \right ) .
}
As before, it follows that $L_{Y^*}(x) > L_{\widehat{Y}}(x)$ contradicting the minimality of $Y^*$.

\pbk 
Case 3: This is similar to Case 2, and hence omitted.
}

Note that a particular consequence of this lemma is that, as claimed, the function $B(M,N,x)$ is a polynomial on each subinterval $[x_m,x_{m+1}]$.

\subsection{A Remez-type algorithm for computing $B(M,N,x)$}
Lemma \ref{l:schonhage} not only gives an expression for $B(M,N,x)$, its proof also suggests a numerical procedure for its computation.  The algorithm follows the steps of the proof and proceeds roughly as follows.  First, a set $Y$ of the form described in Lemma \ref{l:schonhage} is chosen and the polynomial $p = p_{Y,n}$ of Lemma \ref{l:LagrangeEquiosc} is computed.  If \R{schonhage_claim} holds, then, as shown in the proof of Lemma \ref{l:schonhage}, $p(x) = B(M,N,x)$.  If not, then a point $x^* \in \{ x_m \}^{M}_{m=0}$ which maximizes $|p(x_m)|$ is found, and, following the proof once more, a suitable element $y_k$ of $Y$ is exchanged with $x^*$ to construct a new set $\widehat{Y}$.  This process is repeated until \R{schonhage_claim} holds.

\begin{algorithm}[First Remez-type algorithm]
\label{alg:Remez1}
~
\begin{enumerate}
\item Pick a subset $Y \subseteq \{ x_m \}^{M}_{m=0}$ with $|Y| = N+1$ and $x_p,x_{p+1} \in Y$.
\item Compute the polynomial $p = p_{Y,n} \in \bbP_N$ satisfying \R{LagrangeInterpPoly}, where $n$ is such that $x_n = x_p$.
\item Find a point $x^* \in \{ x_m \}^{M}_{m=0}$ with
\bes{
|p(x^*) | = \max_{m=0,\ldots,M} | p(x_m) |.
}
\item If $|p(x^*) | = 1$, then set $B(M,N,x) = p(x)$ and stop.
\item If $|p(x^*) | > 1$, then proceed as follows:
\begin{enumerate}
\item[a)] Suppose that $y_k < x^* < y_{k+1}$ for some $k=0,\ldots,N-1$.  If $\sgn(p(x^*)) = \sgn(p(y_{k}))$ then replace $y_k$ with $x^*$ in the set $Y$, otherwise replace $y_{k+1}$ with $x^*$.
\item[b)] Suppose that $x^* < y_0$.  If $\sgn ( p(x^*)) = \sgn(p(y_0))$ then replace $y_0$ with $x^*$ in the set $Y$, otherwise replace $y_N$ with $x^*$.
\item[c)] Suppose that $x^* > y_N$.  If $\sgn (p(x^*)) = \sgn(p(y_N))$ then replace $y_N$ with $x^*$ in the set $Y$, otherwise replace $y_0$ with $x^*$.
\end{enumerate}
\item Return to step 2.
\end{enumerate}
\end{algorithm}
This algorithm is guaranteed to converge in a finite number of steps.  As shown in the proof of Lemma \ref{l:schonhage}, the exchange performed in step 5 strictly decreases the value $L_Y(x)$.  Since there are only finitely-many possible sets $Y$, the algorithm must therefore terminate in finite time.

In practice, it is usually preferable to exchange more than one point $x^*$ at a time.  This leads to the following algorithm:

\begin{algorithm}[Second Remez-type algorithm]
~
\\
The algorithm is similar to Algorithm \ref{alg:Remez1}, except that in step 3 we find all extrema of $p(x)$ on the set $\{ x_m \}^{M}_{m=0}$.  Note that there are at least $N-3$ such points and at most $N-2$.\footnote{The polynomial $p$ of degree $N$ has a full set of $N$ zeros, hence $N-1$ extrema.  Of those, one extremum is between $x_m$ and $x_{m+1}$, and at most one extremum might be outside the interval under consideration. Hence, the number of extrema that may appear in the algorithm 
is $N-2$ at most, and $N-3$ at least.}
Each subinterval $[y_k,y_{k+1}]$ contains at most one of these extrema.  Hence we now proceed with the exchange as in step 5 above for each such point.
\end{algorithm}

\subsection{Numerical results}
\label{s:Numerical}

\begin{figure}
\begin{center}
{\includegraphics[height=.31\textwidth]{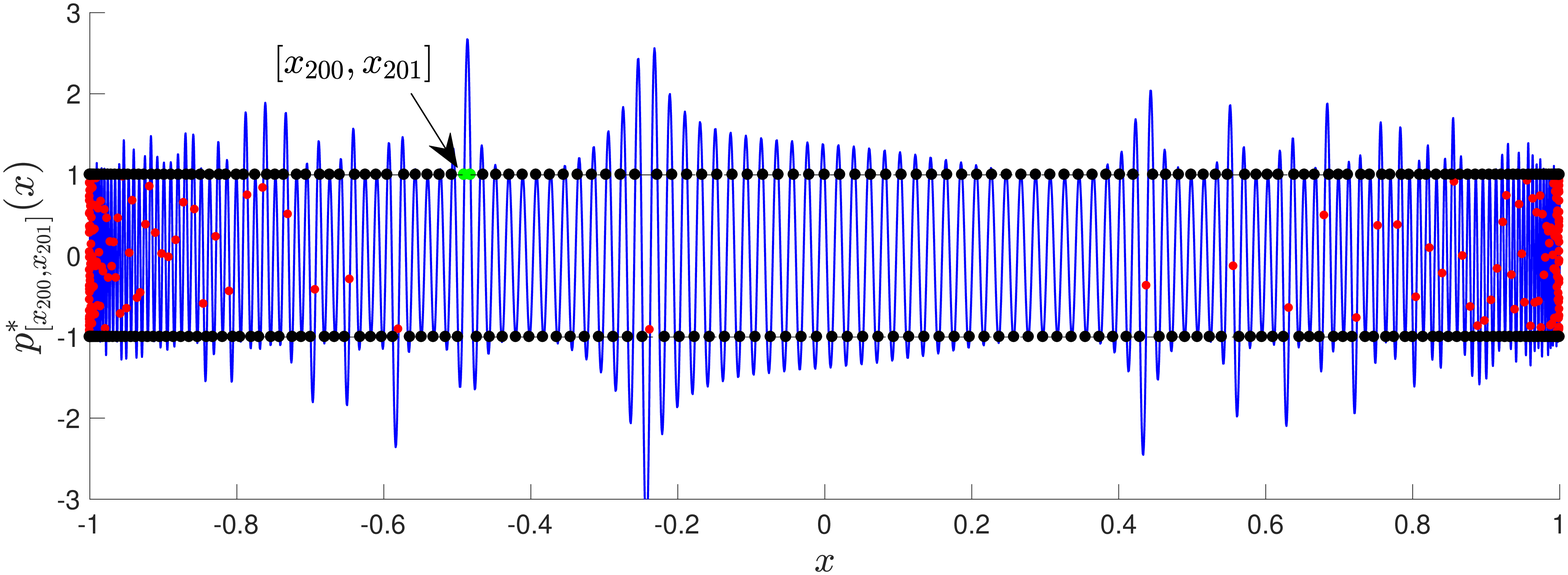}} \\ 
{\includegraphics[height=.25\textwidth]{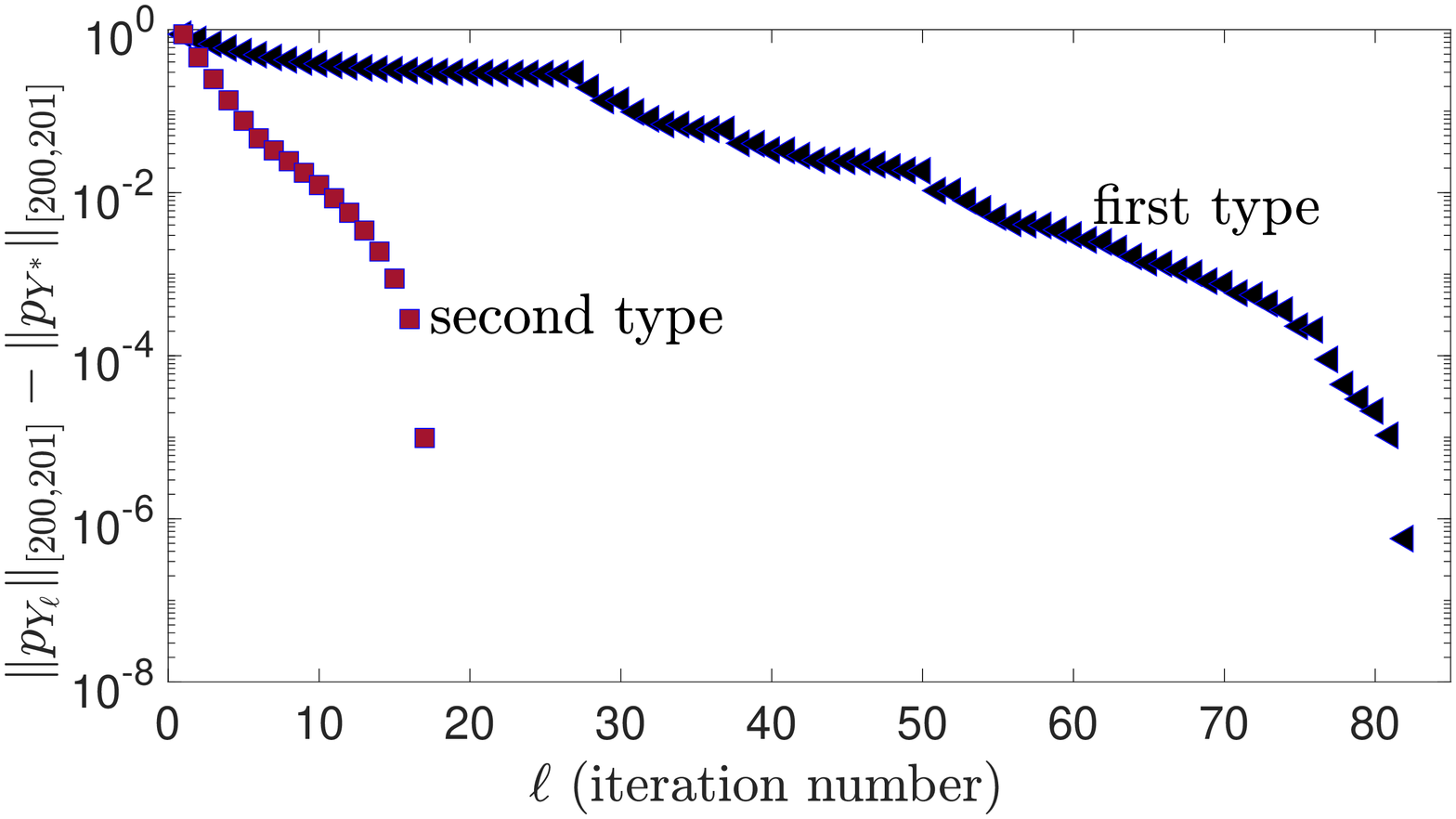}}
\end{center}
\caption{Top: maximal polynomial of degree 300 over the interval $[x_{200},x_{201}]$, where $\{x_m\}_{m=0}^M$ is the set of 500 points drawn from the (OC) density. Dot markers show the value of the polynomial at the grid $x_m$, with black dots corresponding to points for which the polynomial evaluates to $\pm 1$. Bottom: Convergence plot of the first and second Remez-type algorithms to this polynomial.}
\label{fig:RemezExample}
\end{figure}

The performance of the first and second-type Remez algorithms is presented in Figure~\ref{fig:RemezExample}, where the maximal polynomial of degree 300 over the interval $[x_{200},x_{201}]$ is computed for the (OC) case. The first-type algorithm takes over 80 iterations to converge, while the second-type computes the maximal polynomial in 18 iterations. In this experiment, the algorithm was started using the subset $Y \subseteq \{ x_m \}^{500}_{m=0}$ of 301 points closest to the Chebyshev points of the second kind, that is, the mock-Chebyshev subset. The iteration count can be significantly higher when the initial subset of points is selected at random. We also point out that the algorithm may fail to converge if at any iteration the interpolation set $Y$ is very ill-conditioned. In double-precision, the Lebesgue constant for $Y$ must not exceed $10^{16}$. Choosing mock-Chebyshev points to initialize the procedure, therefore, reduces the likelihood of failed iterations. 

\begin{figure}
\begin{center}
\begin{tabular}{c c c}
& $N=10$, $M=20$ & $N=20$, $M=40$ \\
\rotatebox{90}{\hspace{1.2cm}(C2)}&\includegraphics[height=.19\textwidth]{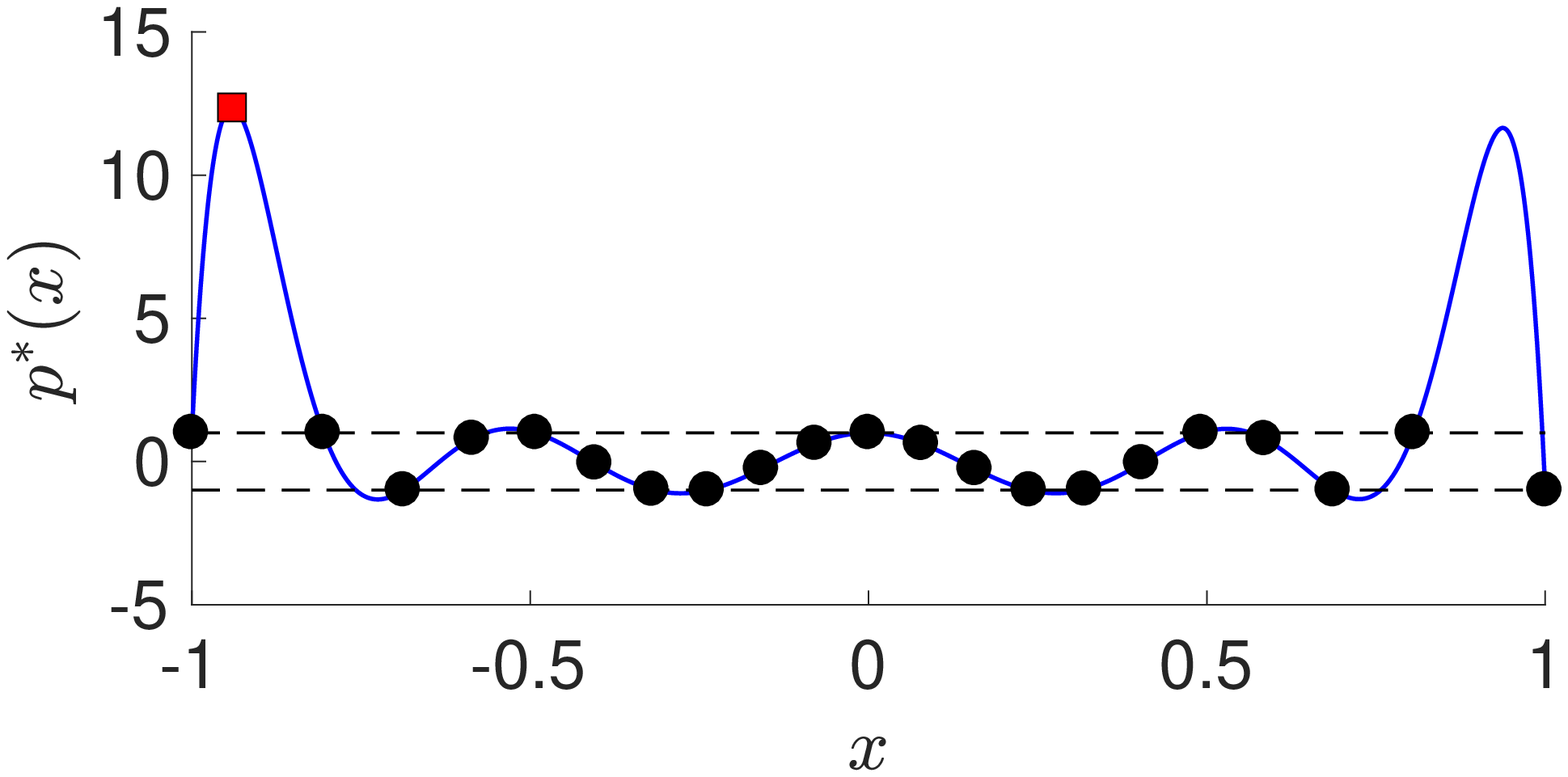} & \includegraphics[height=.19  \textwidth]{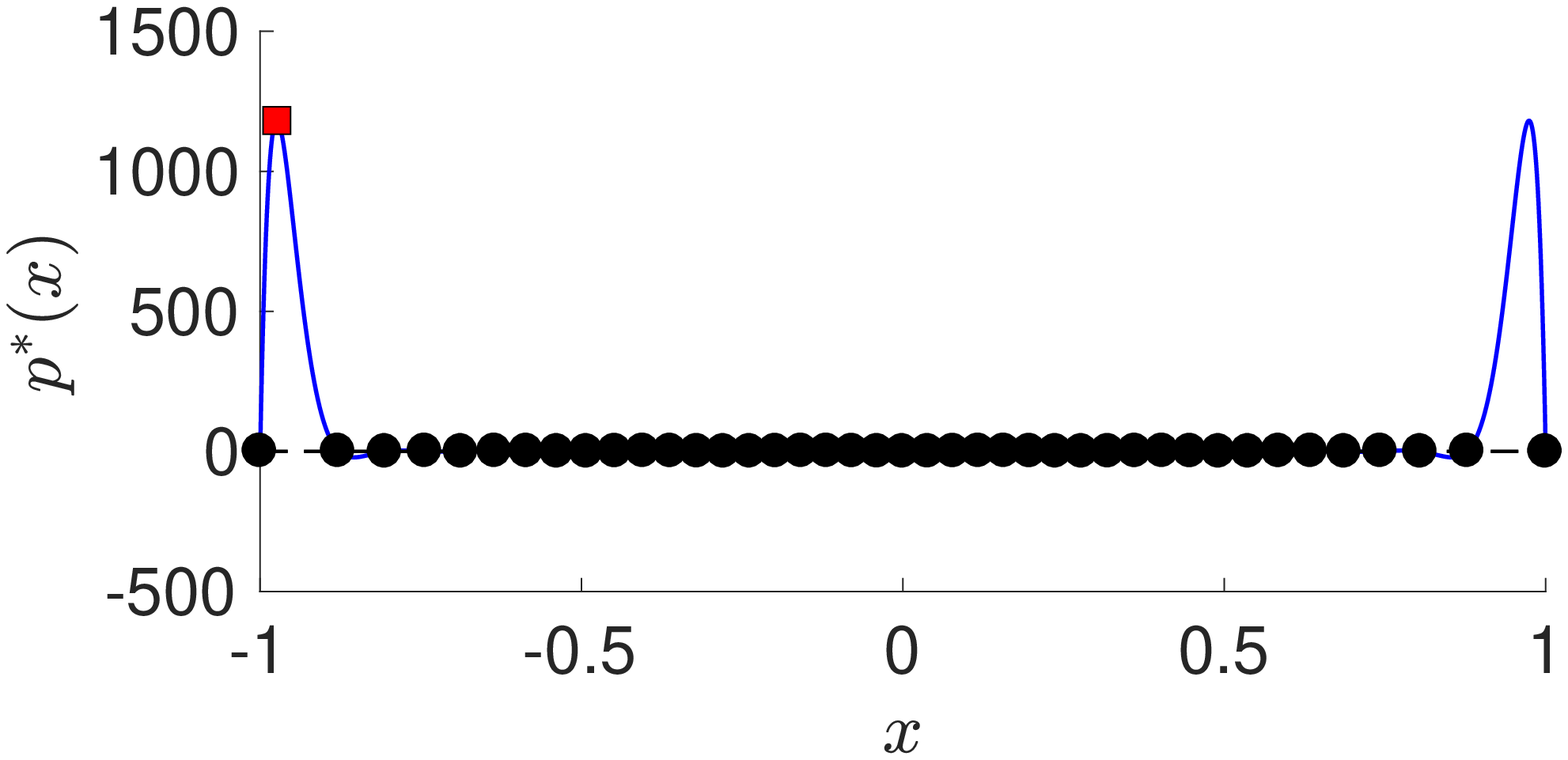}  \\
\rotatebox{90}{\hspace{1.3cm}(U)} &\includegraphics[height=.19\textwidth]{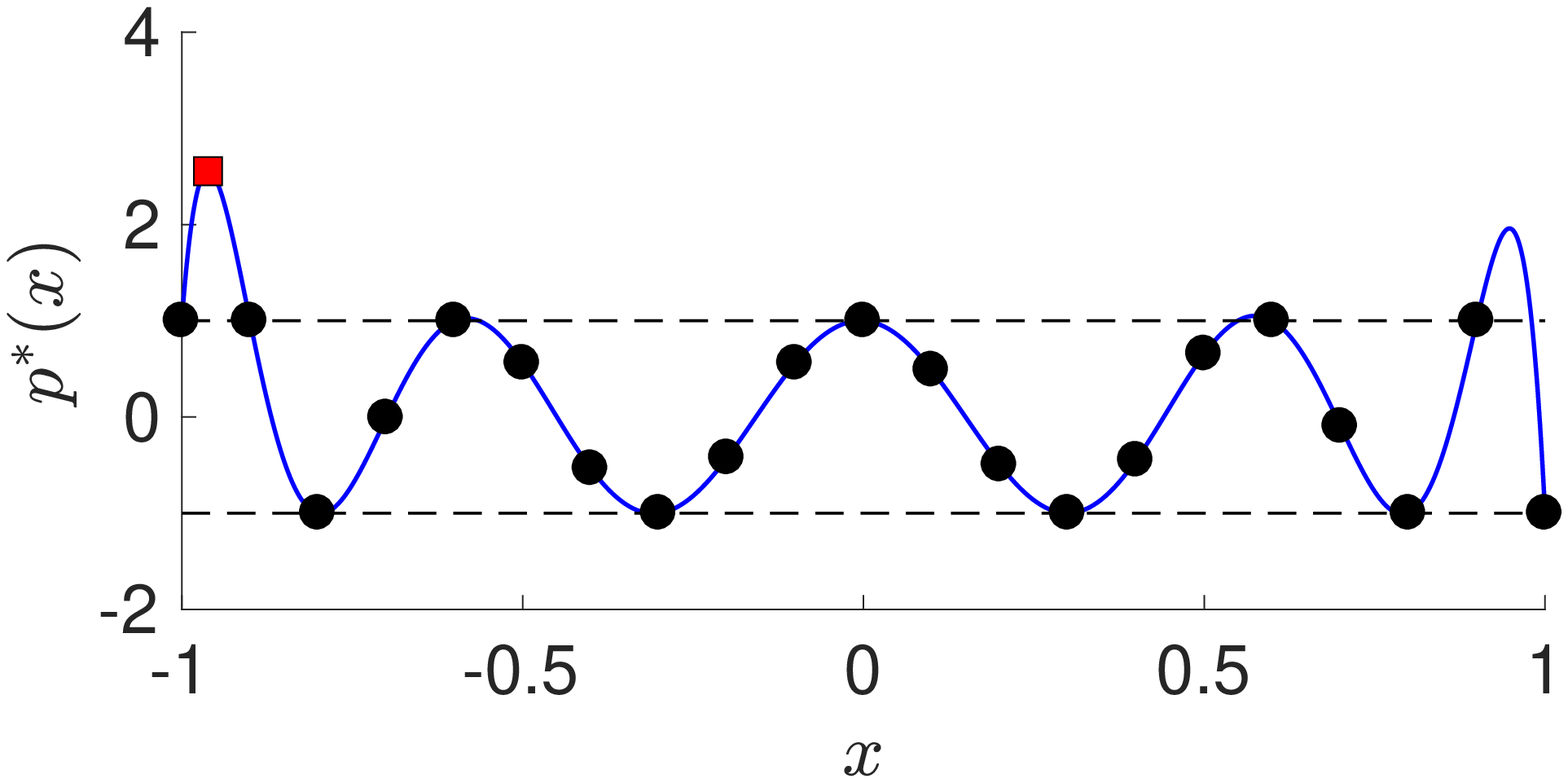} & \includegraphics[height=.19\textwidth]{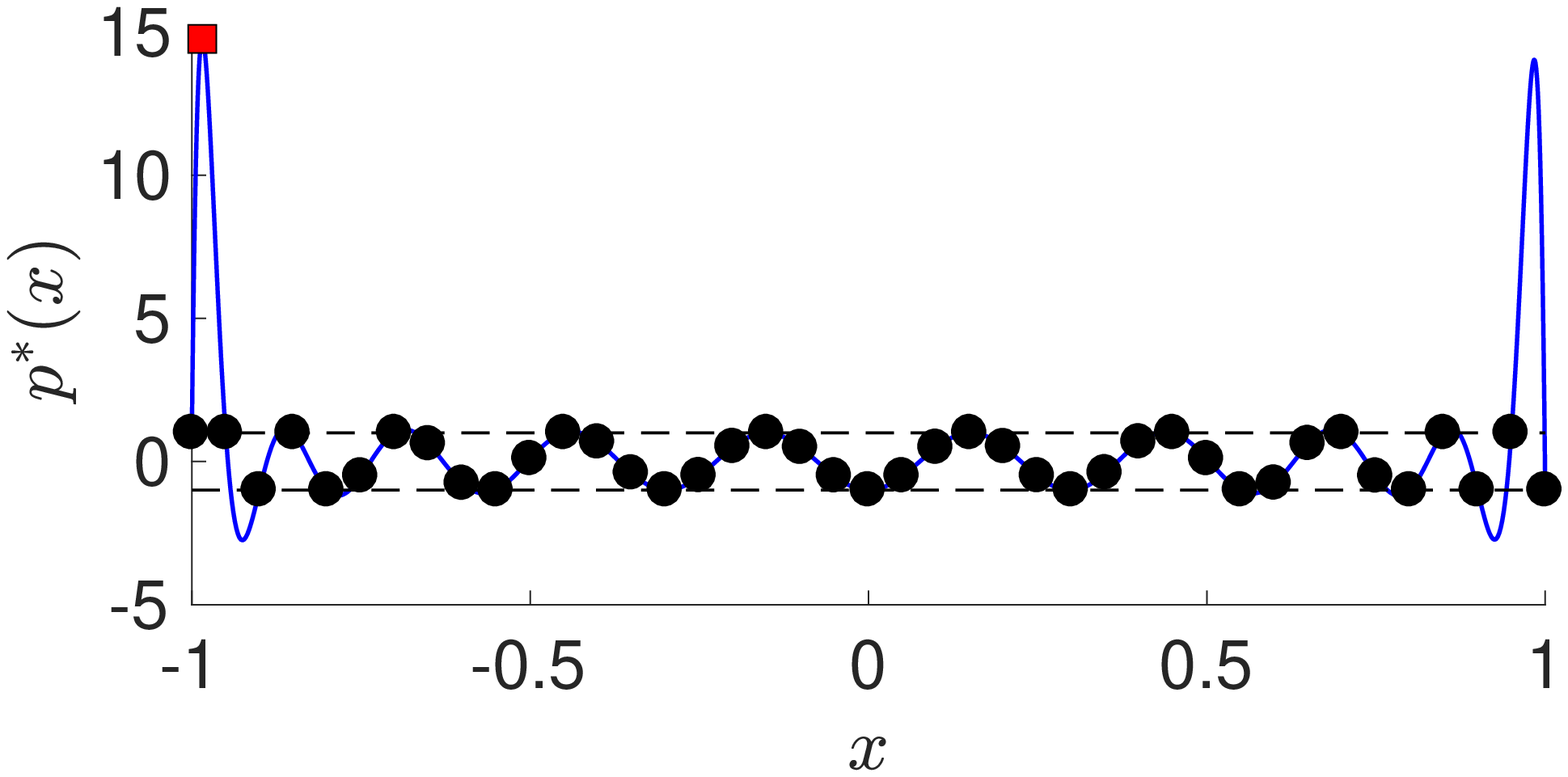}  \\
\rotatebox{90}{\hspace{1.2cm}(UC)}  &\includegraphics[height=.19\textwidth]{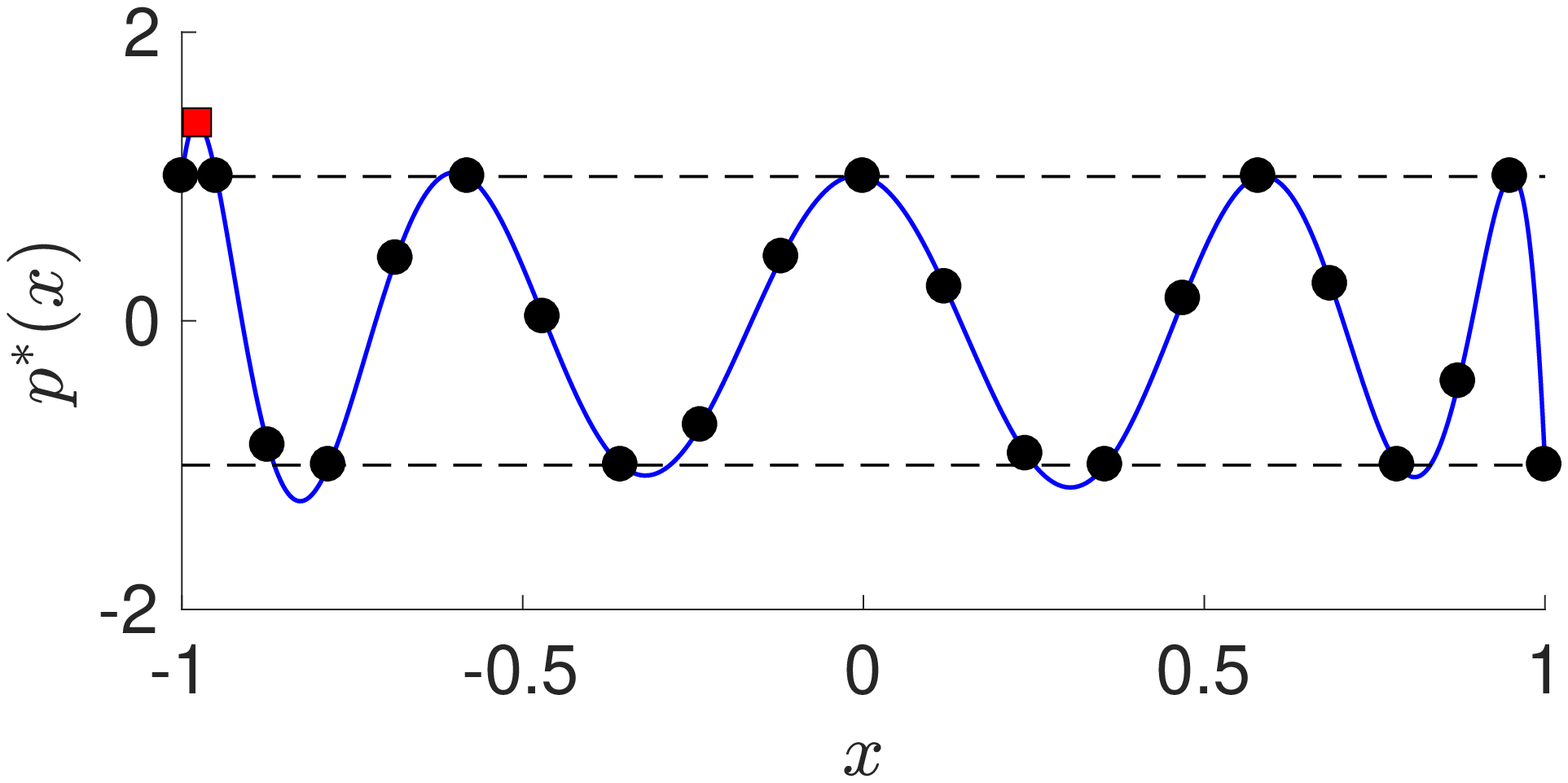} & \includegraphics[height=.19\textwidth]{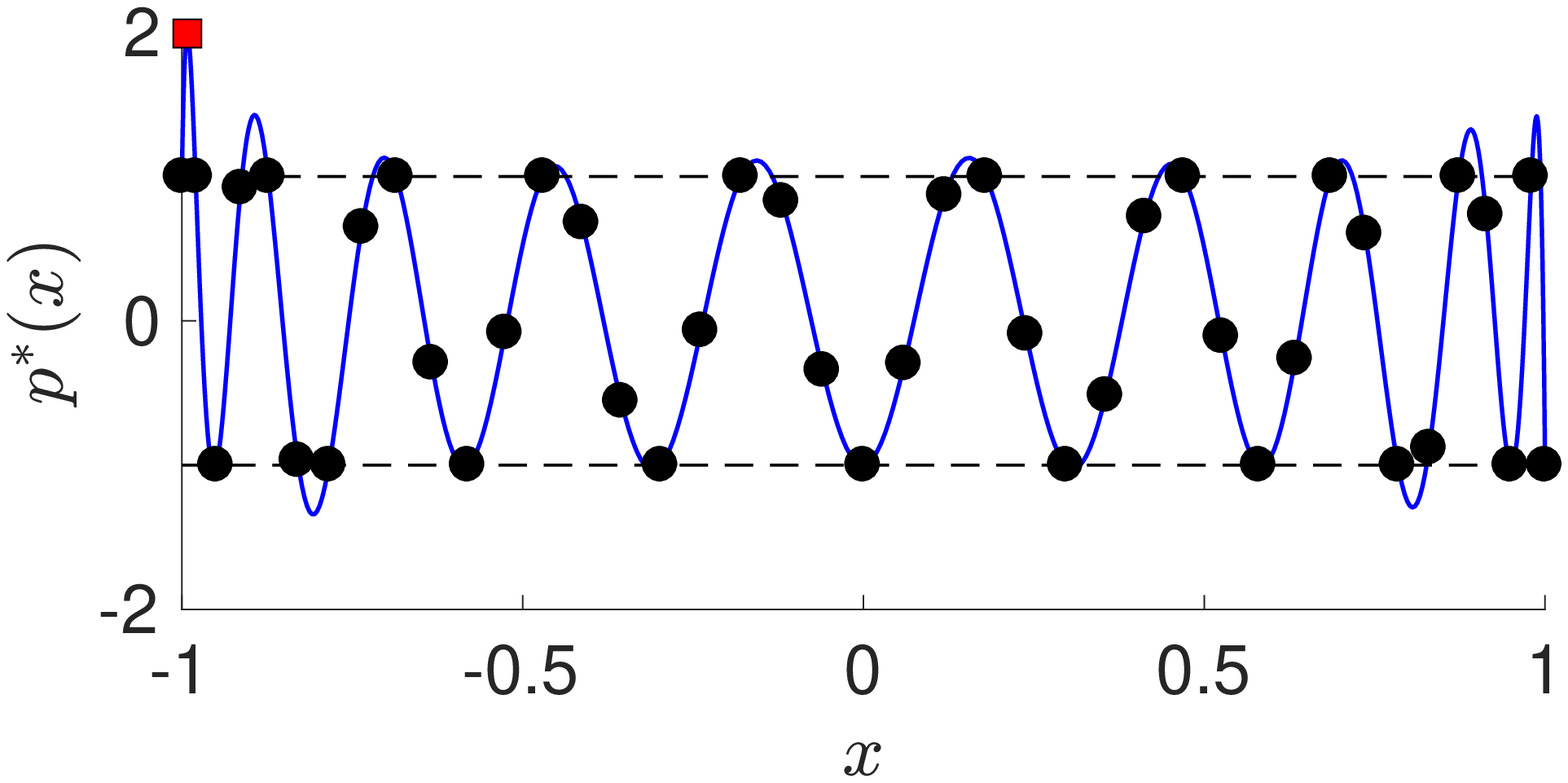}  \\
\rotatebox{90}{\hspace{1.2cm}(OC)}  &\includegraphics[height=.19\textwidth]{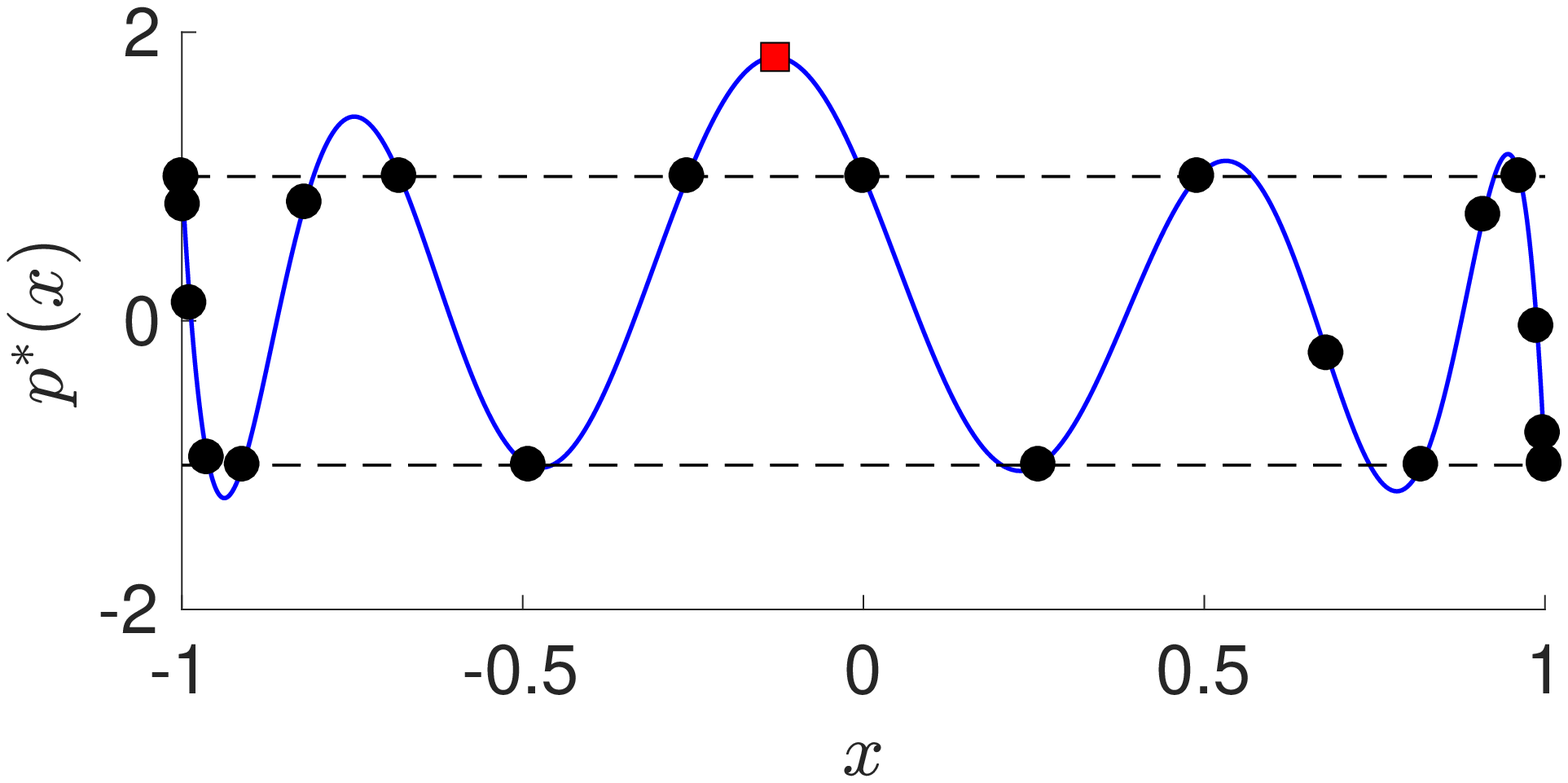} & \includegraphics[height=.19\textwidth]{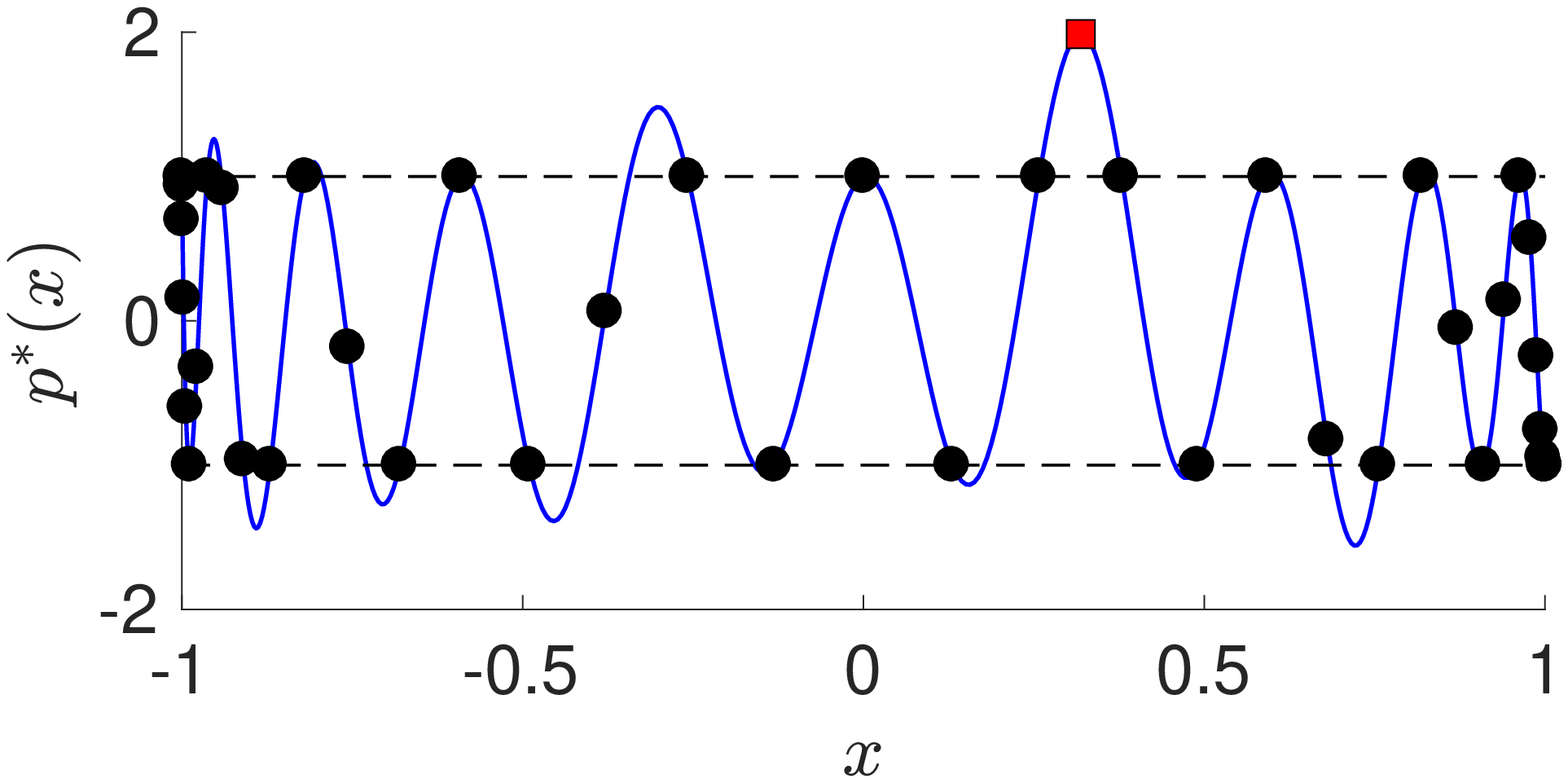}  \\
\end{tabular}
\end{center}
\caption{Maximal polynomials for different node sets. Square markers show the maximum point for each case, while the dot markers show the value of the polynomials at the grid points.}
\label{f:example}
\end{figure}

To find the maximal polynomial over the whole interval, the Remez procedure must be repeated for every subinterval $[x_m, x_{m+1}]$, unless the location where the maximum is achieved is known. In the case of equispaced nodes the maximum is known to near the endpoints.  In the (OC) case, the maximum is in the interior of the interval as illustrated in Figure~\ref{f:example}, but not necessarily at the subinterval closest to the center as shown in the bottom right panel. Although only moderate values of $M$ and $N$ were used in this figure, the Remez algorithm is able to compute maximal polynomials of much larger degrees, as shown in Figures~\ref{fig:prop5_3} and \ref{fig:prop5_3B}.


\section{Concluding remarks}

We have presented a generalized impossibility theorem for approximating analytic functions from $M+1$ nonequispaced points.  This follows from a new lower bound for the maximal behaviour of a polynomial of degree $N$ that is bounded on an arbitrary set of points.  By specializing to modified Jacobi weight functions, we have derived explicit relationships between the parameter $\gamma = \max \{ \alpha , \beta , -1/2 \}$, the rate of exponential convergence and the rate of exponential ill-conditioning.  Polynomial least-squares using a polynomial of degree $N$ transpires to be an optimal stable method in view of this theorem.  In particular, the sampling rate $M \asymp N^{2(\gamma+1)}$, where $\gamma = \max \{ \alpha , \beta , -1/2 \}$ is both sufficient and necessary for stable approximation with optimal convergence.

There are a number of directions for future investigation.  First, we have only derived an upper bound for $B(M,N)$ in the case of ultraspherical weight functions (see Remark \ref{r:Rakhmanov_upper}).  We expect the techniques of \cite{RakhmanovPolyBds} can be extended to the modified Jacobi case whenever $\max \{\alpha,\beta \} > -1/2$.  Second, we have observed numerically that there is an exponential blow-up for $B(M,N)$ in the case $-1 < \gamma < -1/2$ when $M = c N$ for some $c>0$ below a critical threshold.  This remains to be proven.  Third, we have mentioned in passing recent results on sufficient sampling rates when drawing random points from measures associated with modified Jacobi weight functions.  It would be interesting to see if the techniques used in the paper could also establish the necessity (in probability) of those rates.  Fourth, the extension of our results to analytic functions of two or more variables remains an open problem.  Fifth, we note in passing that there is a related impossibility theorem for approximating analytic functions from their Fourier coefficients \cite{AdcockHansenShadrinStabilityFourier} (this is in some senses analogous to the case of equispaced samples).  The extension to nonharmonic Fourier samples may now be possible using the techniques of this paper.  Note that necessary and sufficient sampling conditions for this problem have been proven in \cite{AdcockGataricRomeroNonuniform} and \cite{BAMGACHNonuniform1D,AdcockGataricHansenICOSAHOM} respectively.

Finally, we remark the following.  The impossibility theorems proved here and originally in \cite{TrefPlatteIllCond} assert the lack of existence of stable numerical methods with rapid convergence.  They say nothing about methods for which the error decays only down to some finite tolerance.  If the tolerance can be set on the order of machine epsilon, the limitations of such methods in relation to methods which have theoretical convergence to zero may be of little consequence  in finite precision calculations.  Several methods with this behaviour have been developed in previous works \cite{FEStability,AdcockPlatteMapped}.  The existence (or lack thereof) of impossibility theorems in this finite setting is an open question.

\section*{Acknowledgements}
BA acknowledges the support of the Alfred P.\ Sloan Foundation and the Natural Sciences and Engineering Research Council of Canada through grant 611675.  RBP was supported by NSF-DMS 1522639, NSF-DMS 1502640 and AFOSR FA9550-15-1-0152. We have benefited from using \texttt{Chebfun} (\url{www.chebfun.org}) in the implementation of our algorithms.

\small
\bibliographystyle{abbrv}
\bibliography{SampCondsRefs}

\end{document}